\documentclass[a4paper,12pt,reqno]{amsart}
\usepackage{amssymb,amsthm}
\usepackage{amsmath}
\usepackage{mathtools}
\usepackage{ifthen}
\usepackage{graphicx}
\usepackage{xcolor}
\usepackage{float}
\usepackage{tcolorbox}
\usepackage{caption}
\usepackage{subcaption}
\usepackage{diagbox}
\usepackage{hyperref} % To create clickable links in the PDF
\theoremstyle{plain}

\newtheorem{thm}{Theorem}[section]
\newtheorem*{thm*}{Theorem}

\numberwithin{equation}{section}

\numberwithin{subcase}{case}
\newtheorem{cor}{Corollary}[section]
\newtheorem{lem}{Lemma}[section]
\newtheorem{prop}{Proposition}[section]

\newtheorem{example}{Example}[section]
\newtheorem{defn}{Definition}[section]

\theoremstyle{definition}
\newcounter {own}
\def\theown {\thesection  .\arabic{own}}

\newenvironment{pf}[1][]{%
 \vskip 3mm
 \noindent
 \ifthenelse{\equal{#1}{}}%
  {{\slshape Proof. }}%
  {{\slshape #1.} }%
 }%
{\qed\bigskip}

\DeclareMathOperator{\lcm}{lcm}

\DeclareMathOperator{\sinc}{sinc}
\newcounter{alphabet}

\newcounter{minutes}
\divide\time by 60
\newcounter{hours}
\multiply\time by 60 \addtocounter{minutes}{-\time}
\begin{document}
\bibliographystyle{amsplain}
\title{Construction  of irregular complete interpolation sets for shift-invariant spaces }

%=========================================================================
% \thanks{%$^\dagger$
% File:~AntonyPriya2.tex,
%           printed: 2022-02-20,
%           \thehours.\ifnum\theminutes<10{0}\fi\theminutes}
%=========================================================================
\author{Kumari Priyanka}
\author{A. Antony Selvan$^\dagger$}
\address{Kumari Priyanka, Indian Institute of Technology Dhanbad, Dhanbad-826 004, India.}
\email{priyankak4193@gmail.com}
\address{A. Antony Selvan, Indian Institute of Technology Dhanbad, Dhanbad-826 004, India.}
\email{antonyaans@gmail.com}
\keywords{ B-splines, Complete interpolation sets, Exponential splines, Toeplitz operators, Transversal sets, Shift-invariant spaces.\\
$^\dagger$ {\tt Corresponding author}
}
\maketitle
\pagestyle{myheadings}
\markboth{Kumari Priyanka and A. Antony Selvan}{Construction of irregular complete interpolation..}
\let\thefootnote\relax\footnotetext{An earlier version of this paper was posted on \href{https://arxiv.org/abs/2408.09099}{arXiv:2408.09099}.}
% \footnotetext{This paper was previously posted on \href{https://arxiv.org/abs/2408.09099}{arXiv:2408.09099}.}
\begin{abstract}
For several shift-invariant spaces, there exists a real number $a\in\mathbb{R}$ such that the set $a+\mathbb{Z}$ is a complete interpolation set. In this paper, we characterize the complete interpolation property of 
the set $(a+\mathbb{N}_0)\cup(\alpha+a+\mathbb{N}^{-})$ for shift-invariant spaces using  Toeplitz operators.  Using this
characterization, we  determine all $\alpha$ for which the sample set $\mathbb{N}_0\cup\alpha+\mathbb{N}^{-}$ forms
a complete interpolation set for transversal-invariant spaces.
We introduce a new recurrence relation for exponential splines, examines the zeros of these splines, and explores the zero-free region of the doubly infinite Lerch zeta function. Consequently, we demonstrate that $\left\langle\frac{m}{2}\right\rangle+\mathbb{N}_0\cup\alpha+\left\langle\frac{m}{2}\right\rangle+\mathbb{N}^{-}$ is a complete interpolation set for a  shift-invariant spline space of order $m\geq 2$ if and only if $|\alpha|<1/2$.
\end{abstract}

\section{Introduction}
Sampling and interpolation in shift-invariant spaces have been a central topic in mathematical analysis, particularly due to their applications in signal processing, approximation theory, and functional analysis. A classical example of a shift-invariant space is the space of bandlimited signals with bandwidth $\pi.$ 
In the literature, an almost characterization of sampling sets and of  interpolation sets for a large class of shift-invariant spaces has been given in terms of Beurling densities
\cite{Beurling, Grochenig, Kahane, Landau}. The problem of complete interpolation sets  is much more subtle. Pavlov \cite{Pavlov} characterized complete interpolation sets for the space of bandlimited functions in the context of non-harmonic Fourier series.  Recently, the authors in \cite{Baranov} solved the problem of  complete interpolating sets for the Gaussian shift-invariant space.

 Characterizing complete interpolation sets for the space of multi-bandlimited signals is a challenging problem and constructing such a set is generally quite difficult.
 This problem is equivalent to finding Riesz bases of exponential systems on a disconnected set $\Omega$. This problem is solved for a finite union of disjoint intervals in \cite{Kozma} and bounded remainder sets in \cite{Grepstad}. We refer to \cite{Olevskii} for an extensive survey of recent progress on the exponential bases problem.

In many concrete shift-invariant spaces, it has been established that there exists a real number $a\in\mathbb{R}$ such that the set $a+\mathbb{Z}$ forms a complete interpolation set. See \cite{AntoRad, Torres, Janssen, Schoenberg, Walter1}. The idea of shifting lattices traces back to Kohlenberg's work \cite{Kohlenberg}, where he first considered the conditions under which a union of shifted lattices could form a sampling and interpolation set for Paley-Wiener spaces  $\mathcal{PW}(\Omega),$ especially when $\Omega$ is a union of two intervals of equal length.  This foundational work was further extended by Lyubarskii and Seip \cite{Lyubarskii}. Schoenberg \cite{Schoenberg} proved that $\mathbb{Z}$ is a complete interpolation set for shift-invariant spline spaces of even order, whereas $\frac{1}{2}+\mathbb{Z}$ is a complete interpolation set for shift-invariant spline spaces of odd order.  In \cite{Massopust}, the authors investigated the problem of complete interpolation  for complex B-splines. 

In this paper, we are interested in determining all $\alpha$ for which the sample set $(a+\mathbb{N}_0)\cup(\alpha+a+\mathbb{N}^{-})$ forms a complete interpolation set for shift-invariant spaces, when $a+\mathbb{Z}$ is a complete interpolation set.  This problem has connections to the study of Riesz bases, particularly in the context of windowed exponentials. Lai \cite{Lai} explored whether the set of complex exponential $E(\Lambda)=\{e^{2\pi i\lambda x}:\lambda\in\mathbb{N}_0\cup\alpha+\mathbb{N}^{-}\}$ forms a Riesz basis for $L^2[0,1].$  He proved that $E(\Lambda)$ is a Riesz basis for $L^2[0,1]$  if and only if $|\alpha|<1/2$. This problem can also be viewed from the perspective of perturbation theory. While there have been significant advances in understanding perturbations of the standard exponential basis on intervals, our results extend beyond the current literature, particularly concerning exponential bases on transversal sets and arbitrary shifts of the negative frequencies. In this paper, we completely characterize all $\alpha$ for which the sample set $(a+\mathbb{N}_0)\cup(\alpha+a+\mathbb{N}^{-})$ forms
a complete interpolation set for transversal-invariant and spline spaces. It is known that  $\mathbb{Z}$ is  a complete interpolation set  for transversal-invariant spaces (see \cite{Torres}). Our result reveals that the set of all 
$\alpha$ for which $\mathbb{N}_0\cup\alpha+\mathbb{N}^{-}$ 
forms a complete interpolation set for a transversal-invariant space is more complex than in the classical case $E=[0,1]$.
Although our problem pertains to sampling theory,
it leads us to investigate the zeros of Schoenberg's exponential splines and the zero-free region of the doubly infinite Lerch zeta function.

In the study of polynomial sequences, recurrence relations play a vital role for studying the locations of zeros. Schoenberg provided a recurrence relation for exponential splines in \cite{Schoenberg1}. However, his recurrence relation does not assist us in solving our problem.  This paper introduces a new recurrence relation for exponential splines, which provides new insights into their zeros and will be used in constructing complete interpolation sets.

The paper is organized as follows: Section 2 provides key definitions, relevant results, an overview of Toeplitz operators along with related theorems that will be utilized in rest of the paper.
In Section 3, we completely characterize all $\alpha$ for which the sample set $\mathbb{N}_0\cup\alpha+\mathbb{N}^{-}$ forms
a complete interpolation set for transversal invariant spaces (Theorem \ref{general}).
 Section 4 introduces a new recurrence relation for exponential splines (Theorem \ref{phigm}), examines the zeros of these splines (Corollary \ref{zeroofphi}), and explores the zero-free region of the doubly infinite Lerch zeta function (Corollary \ref{zerooflerch}). Finally, in  Section 5, we prove that $\left\langle\frac{m}{2}\right\rangle+\mathbb{N}_0\cup\alpha+\left\langle\frac{m}{2}\right\rangle+\mathbb{N}^{-}$ is a complete interpolation set for a  shift-invariant spline space of order $m\geq 2$ if and only if $|\alpha|<1/2$ (Theorem \ref{cisspline}).
\section{Preliminary Results}
A sequence of vectors $\{f_n:~n\in\mathbb{Z}\}$ in a separable Hilbert space $\mathcal{H}$ is said to be a Riesz sequence if there exist two positive constants $A$ and $B$ such that
\begin{equation}
A\sum_{n\in\mathbb{Z}}|d_n|^2\leq\big\|\sum_{n\in\mathbb{Z}}d_nf_n\big\|^2_{\mathcal{H}}\leq
B\sum_{n\in\mathbb{Z}}|d_n|^2,
\end{equation}
for every $(d_n)\in\ell^2(\mathbb{Z})$;
it is said to be a {\it frame} if there
exist constants $0< A \leq B<\infty$ such that
\begin{equation}
A\|f\|_\mathcal{H}^2\leq\displaystyle\sum_{n\in\mathbb{Z}}|\langle
f,f_n\rangle_\mathcal{H}|^2\leq B\|f\|_\mathcal{H}^2,
\end{equation}
for every $f\in\mathcal{H}$.
A sequence is called a Riesz basis if it is a frame and Riesz sequence.
The shift-invariant space $V(\phi)$ is defined as 
\begin{equation*}
V(\phi):=\left\{\sum\limits_{k\in\mathbb{Z}}c_k\phi(\cdot-k):(c_k)\in\ell^2(\mathbb{Z})\right\}\subset L^2(\mathbb{R}).
\end{equation*}
The generator $\phi$ is said to be stable if $\{\phi(\cdot-n) : n \in \mathbb{Z}\}$ is a Riesz basis for $V(\phi)$
In addition, if $\phi$ is a continuous function such that
for some 
$\epsilon>0$, $$\phi(x)=\mathcal{O}(|x|^{-0.5-\epsilon})~ \text{as}~ x\to\pm\infty,$$
then we say that $\phi$ is a stable generator for $V(\phi)$.
In this case, every $f\in V(\phi)$  can be written as 
\begin{eqnarray*}
f(x)=\sum\limits_{k\in\mathbb{Z}}c_k\phi(x-k).
\end{eqnarray*}
It is well known that $\phi$ is stable  if and only if there exist two positive constants $A, B$ such that
\begin{equation}\label{eqn2.5}
 A\leq\displaystyle\sum_{n\in\mathbb{Z}}|\widehat{\phi}(\xi+n)|^2\leq B~~ a.e.~~\xi\in\mathbb{R},
\end{equation}
where $\widehat{\phi}$ is the Fourier transform of $\phi$, and is
given by
$$\widehat{\phi}(w)= \int_{-\infty}^{\infty} \phi(t) e^{-2\pi iwt}dt,~ w\in\mathbb{R}.$$
\begin{defn}
Let $\phi$ be a stable generator for $V(\phi).$ Then a set $\Gamma= \{x_n: n\in\mathbb{Z}\}$ of real numbers is said to be 
\begin{enumerate}
\item[$(i)$] a set of stable sampling for $V(\phi)$ if there exist constants $A, B>0$ such that 
\begin{equation*}
A\|f\|^2_2\leq\sum\limits_{n\in\mathbb{Z}}|f(x_n)|^2\leq B\|f\|_2^2,~
\text{for all}~f\in V(\phi);
\end{equation*}
\item[$(ii)$] a set of interpolation for $V(\phi)$ if the interpolation problem $$f(x_n)=c_{n},~n\in\mathbb{Z},$$ has a solution $f\in V(\phi)$ for every square summable sequence $\{c_{n}: n\in\mathbb{Z}\}.$
\end{enumerate}
If $\Gamma$ is both set of stable and interpolation for $V(\phi),$ then it is called a complete interpolation set for $V(\phi).$
\end{defn}
These notions can be expressed in terms of the pre-Grammian matrix $$U:=[\phi(x_n-k)]_{x_n\in\Gamma, k\in\mathbb{Z}}$$
of the generator $\phi$ and the sample set $\Gamma$. See \cite{Grochenig}.
\begin{prop}\label{prop1.1}
Let $\phi$ be a stable generator for $V(\phi)$ and $\Gamma= \{x_n: n\in\mathbb{Z}\}\subset\mathbb{R}$. Then
\begin{itemize}
\item [$(i)$] $\Gamma$ is a set of stable sampling for $V(\phi)$  if and only if  $U$ is bounded above and below on $\ell^2(\mathbb{Z})$.
\item[$(ii)$]  $\Gamma$ is a set of interpolation for $V(\phi)$ if and only if $U$ is a bounded surjective operator from $\ell^2(\mathbb{Z})$  onto $\ell^2(\mathbb{Z})$.
\item[($iii)$] $\Gamma$ is a complete interpolation set for $V(\phi)$ if and only if $U$ is  bounded invertible on $\ell^2(\mathbb{Z})$.  
\end{itemize}
\end{prop}

The Zak transform $\mathcal{Z} f$ of a function $f\in L^2(\mathbb{R})$ is defined by 
\begin{equation}\label{zakdef}
\mathcal{Z}f(x,y):=\sum\limits_{k\in\mathbb{Z}}f(x- k)e^{2\pi iky},
~~a.e.,~x,y\in\mathbb{R}.
\end{equation}
The Zak transform of a function $f$ and its Fourier transform $\widehat{f}$ is related by
\begin{eqnarray}\label{zak}
\mathcal{Z}f(x,y)=e^{2\pi ixy}\mathcal{Z}\widehat{f}(y,-x),~~a.e.,~x,y\in\mathbb{R}.
\end{eqnarray}
The Zak transform is a unitary operator from $L^2(\mathbb{R})$ onto $L^2([0,1]^{2})$. We refer to \cite{Grochenig1} for other properties of the Zak transform.

A set $E\subseteq\mathbb{R}$ 
is called a $1$-transversal if the collection of translates 
$\{E+n:n\in\mathbb{Z}\}$ forms a measurable partition of $\mathbb{R}$.
In other words, the union of these translates covers $\mathbb{R}$
 without overlap, except possibly on a set of measure zero.
It is well known that
$\left\{e^{2\pi in \cdot} :n\in\mathbb{Z}\right\}$ forms an orthonormal basis
for $L^2(E)$.
For each $F\in L^2(E)$, we define the operator $T$ by
\begin{equation}\label{eqn2.11}
(TF)(x):= \sum\limits_{n\in\mathbb{Z}}\left\langle F,e^{2\pi
in\cdot}\right\rangle_{L^2(E)}\phi(x+n)=\left\langle F,\overline{\mathcal{Z}\phi(x,.)}\right\rangle_{L^2(E)}, ~~x\in\mathbb{R}.
\end{equation}
It is clear that $TF\in V(\phi)$. Since $(Te^{2\pi
i m\cdot})(x)=\phi(x+m)$, $T$ is
a bounded invertible linear operator from $L^2(E)$ onto
$V(\phi)$. Since 
$$f(\lambda_n)
=\left\langle T^{-1}f,\overline{\mathcal{Z}\phi(x_n,\cdot)}\right\rangle_{L^2(E)},
$$ 
we can prove the following result, which is a slight refinement of the specific case $E=[0,1]$ analyzed in  \cite{AntoRad, Garcia}. 
\begin{prop}\label{zakthm}
Let $\Gamma= \{x_n: n\in\mathbb{Z}\}\subset\mathbb{R}$
and $\phi$ be a stable generator for $V(\phi)$. 
Then
\begin{itemize}
\item [$(i)$] $\Gamma$  is a set of stable sampling for $V(\phi)$ if and only if  $\{\mathcal{Z}\phi(x_n,\cdot):n\in\mathbb{Z}\}$ is a frame for $L^2(E)$.
\item [$(ii)$]  $\Gamma$  is a set of interpolation for $V(\phi)$ if and only if  $\{\mathcal{Z}\phi(x_n,\cdot):n\in\mathbb{Z}\}$ is a Riesz sequence for $L^2(E)$.
\item [$(iii)$] $\Gamma$ is a complete interpolation set for $V(\phi)$ if and only if  $\{\mathcal{Z}\phi(x_n,\cdot):n\in\mathbb{Z}\}$ is a Riesz basis for $L^2(E)$.
\end{itemize}
\end{prop}

A bounded linear operator $\mathcal{A}$  on $\ell^2(\mathbb{Z})$ associated with a matrix $\mathcal{A}=[a_{rs}]$ is said to be a \textit{Laurent operator} if $a_{r-k,s-k}=a_{r,s}$, for every $r$, $s$, $k$
$\in\mathbb{Z}$.
For $m\in L^\infty[0,1]$, let us define $\mathcal{M}:L^2[0,1]\to
L^2[0,1]$ by \[(\mathcal{M}f)(t):=m(t)f(t), ~\text{for}~a.e.~t\in [0,1].\]
Then the operator 
$\mathcal{A}=\mathcal{F}\mathcal{M}\mathcal{F}^{-1}$  is  called a Laurent operator defined by the symbol $m$, where
$\mathcal{F}$ denotes the
the Fourier transform
from  $L^2(\mathbb{T})$ onto $ \ell^2(\mathbb{Z})$, and is  defined by $$(\mathcal{F}f)(n)=\widehat{f}(n)=\int_{0}^{1} f(t) e^{-2\pi\mathrm{i}n t}~dt,~n\in\mathbb{Z}.$$
In the sequel when there is no confusion possible, we use $\widehat{f}$ to denote either the Fourier coefficient or Fourier transform of a function.

If a Laurent operator $\mathcal{A}$ is invertible, then $\mathcal{A}^{-1}$ is also a Laurent operator with the symbol $\tfrac{1}{m}$ and the matrix of $\mathcal{A}^{-1}$ is given by
$$\mathcal{A}^{-1}=\left[\widehat{\left(\tfrac{1}{m}\right)}(r-s)\right].$$

Let $\mathbb{N}_0=\mathbb{N}\cup\{0\}$ and  $\mathbb{N}^{-}$ be the
set of all negative integers. 
A bounded linear operator $\mathcal{T}:\ell^2(\mathbb{N}_0)\to\ell^2(\mathbb{N}_0)$ is a Toeplitz operator if, with respect to the standard basis of $\ell^2(\mathbb{N}_0),$ $\mathcal{T}$ admits matrix representation 
\begin{eqnarray*}
 \mathcal{T}=\left(
\begin{array}{ccccccc}
 a_0 & a_{-1} & a_{-2} & \cdots \\
 a_1 & a_{0}  & a_{-1} & \cdots \\
 a_2 & a_{1}  & a_{0}  & \cdots \\
\vdots & \vdots  & \vdots & \ddots \\ 
 \end{array}
\right).
\end{eqnarray*}

Let $\ell_2 = \Big\{(a_n) \in \ell^2(\mathbb{Z}):a_n=0 \text{ for all } n \in \mathbb{N}^{-}\Big\}$. 
Consider the orthogonal projection $P$ of $\ell^2(\mathbb{Z})$ onto $\ell_2$ defined by
$$P(\dots,a_{-1}, a_0, a_1,\dots)=(\dots,0,0, a_0, a_1,\dots).$$
For a given Laurent operator $\mathcal{A}$ on $\ell^2(\mathbb{Z})$ defined by $m$, the operator $\mathcal{T}$  defined by 
$\mathcal{T}=P\mathcal{A}|\ell_2$ is a Toeplitz operator. In this case, we say that $\mathcal{T}$ is a Toeplitz operator defined by the symbol $m$. 

The invertibility of a Toeplitz operator is characterized by the properties of its symbol.
To be precise,
let $\gamma:[0, 1]\to\mathbb{C}$ be a closed curve not passing through the point $z_0.$ Then there exists a continuous function $\sigma:[0, 1]\to\mathbb{C}$ such that $\gamma(x)=e^{\sigma(x)}+z_0,$ for all $x\in[0, 1].$ The winding number (or Index) of $\gamma$ relative to $z_0,$ denoted by $\textrm{Ind}(\gamma,z_0),$ is defined by
$$\textrm{Ind}(\gamma,z_0):=\frac{\sigma(1)-\sigma(0)}{2\pi i}.$$
If $\gamma$ is a rectifiable curve, then
$$\textrm{Ind}(\gamma,z_0)=\frac{1}{2\pi i}\int\limits_{\gamma}\frac{1}{z-z_0}~dz.$$
%For further details see \cite{Ahlfors, Gohberg, Sheil}.
\begin{thm}\cite{Gohberg}
A Toeplitz operator
$\mathcal{T}$ defined by a continuous symbol $m$ is invertible if and only if $m(x)\neq 0$ for all $x\in[0, 1]$ and $\textrm{Ind}(m, 0)=0$.
\end{thm}
In case $m$ is a piecewise continuous function, an invertibility criterion for $\mathcal{T}$ can be stated as follows.
\begin{thm}\cite{Bottcher, Gohberg1}\label{toeplitz1}
Let $\mathcal{T}$ be a Toeplitz operator
defined by a piecewise continuous symbol with jumps at $z_1, \dots, z_n$. Then  
$\mathcal{T}$  is  invertible if and only if $m^{\#}(x)\neq0$ for all $x\in[0, 1]$ and  $\textrm{Ind}(m^{\#},0)=0,$ where
$m^{\#}$
is the closed, continuous curve obtained by connecting the values 
$m(z_j-0)$ and $m(z_j+0)$ with line segments at each jump  $z_j$
forming a naturally oriented curve that represents the essential range of $m$.
\end{thm}
\begin{prop}\cite{Lai}\label{toeplitz2}
Let $T:\ell^2(\mathbb{Z})\to\ell^2(\mathbb{Z})$ be a linear operator such that $T$ has a matrix representation of the form 
$$\left(\begin{array}{c|c}
  A   & B \\
  \hline
  \mathbf{0} & \mathcal{I}   
\end{array}\right),$$
where $A:\ell^2(\mathbb{N}^{-})\to\ell^2(\mathbb{N})$ and $B:\ell^2(\mathbb{N}^{-})\to\ell^2(\mathbb{N}_0)$ are bounded linear operators. Then 
\begin{itemize}
\item [$(i)$] $T$ is a bounded linear operator on $\ell^2(\mathbb{Z}).$
\item[$(ii)$] $A$ is normed bounded below on $\ell^2(\mathbb{N}^{-})$ if and only if $T$ is also normed bounded below. 
\item[$(iii)$] $A$ is invertible on $\ell^2(\mathbb{N}^{-})$ if and only if $T$ is invertible.
\end{itemize}
\end{prop}
Assume that $\phi$ is a stable generator for $V(\phi)$. For each fixed $a\in[0, 1),$ consider the sample set $\Gamma_a=\left\{a+n: n\in\mathbb{Z}\right\}.$ If $\Gamma_a$ is a complete interpolation set for $V(\phi),$ then every $f\in V(\phi)$ can be written as
\begin{equation}\label{theta}
f(x)=\sum_{n\in\mathbb{Z}}f(a+ n)\Theta(x-n),~\Theta(x)=\sum\limits_{v\in\mathbb{Z}}\widehat{\left(\frac{1}{\Psi^\dagger}\right)}(v)\phi(x-v),
\end{equation}
where $\Psi^\dagger(x):=\sum\limits_{n\in\mathbb{Z}}\phi(a+n)e^{2\pi inx}\neq0,$ for all $x\in[0, 1].$ 
Furthermore, $\Theta$ satisfies the interpolation property:
\begin{equation}\label{interpolation}
\Theta(a+n)=\delta_{n0},\text{ for all } n\in\mathbb{Z},
\end{equation}
and its Fourier transform is given by
\begin{equation}\label{fttheta}
\widehat{\Theta}(\xi)=\widehat{\phi}(\xi)\Psi^\dagger(-\xi)^{-1}, ~~a.e.,~ \xi\in\mathbb{R}.  
\end{equation}
For a detailed proof of this result, please refer to \cite{AntoRad, Janssen, Walter1}.

For each $\alpha\in\mathbb{R},$ consider the sample set $\Lambda_{\alpha}=\{a+n:n\in\mathbb{N}_0\}\cup\{\alpha+a+n:n\in\mathbb{N}^{-}\}$ for $V(\phi).$
Let us construct an infinite system
\begin{equation}\label{system}
\begin{aligned}
f(\alpha+a+n)
&=\sum\limits_{k=-\infty}^{-1}d_k\Theta(\alpha+a+n-k)+\sum\limits_{k=0}^{\infty}d_k\Theta(\alpha+a+n-k),~n\in\mathbb{N}^{-},\\ 
f(a+n)&=\sum\limits_{k=-\infty}^{-1}d_k\Theta(a+n-k)+\sum\limits_{k=0}^{\infty}d_k\Theta(a+n-k),~~n\in\mathbb{N}_0.
\end{aligned}    
\end{equation}
Let us denote
\begin{eqnarray}
F=\begin{bmatrix}
\vdots\\
f(\alpha+a-2)\\
f(\alpha+a-1)\\
f(a)\\
f(a+1)\\
f(a+2)\\
\vdots
\end{bmatrix},~
D=\begin{bmatrix}
\vdots\\
d_{-2}\\
d_{-1}\\
d_0\\
d_1\\
d_2\\
\vdots
\end{bmatrix}\text{ and define a matrix } 
U=\left(\begin{array}{c c}
  U_{11}   & U_{12} \\
  U_{21}   & U_{22}
\end{array}\right), \nonumber
\end{eqnarray}
whose entries are 
\begin{align*}
U_{11}&=\left[\Theta(\alpha+a+n-k)\right]_{n,k\in\mathbb{N}^{-}},
~U_{12}=\left[\Theta(\alpha+a+n-k)\right]_{n\in\mathbb{N}^{-},k\in\mathbb{N}_0},\\ 
U_{21}&=\left[\Theta(a+n-k)\right]_{n\in\mathbb{N}_0,k\in\mathbb{N}^{-}},\text{ and }U_{22}=\left[\Theta(a+n-k)\right]_{n,k\in\mathbb{N}_0}.
\end{align*} 
Then we can write \eqref{system} as
\begin{equation}\label{ud=f}
UD=F.
\end{equation}
Using the interpolation property of  $\Theta$,  the operator $U$ has the following block structure:
$$U=\left(\begin{array}{c|c}
 U_{11}    & U_{12} \\ 
 \hline
 \mathbf{0}   & \mathcal{I}
\end{array}\right),$$
where $\mathbf{0}$ and $\mathcal{I}$ denote the zero and identity operators respectively and $U_{11}$ is a Toeplitz operator with the symbol 
\begin{eqnarray*}
\Psi_{\alpha}(x)=\sum\limits_{n\in\mathbb{Z}}\Theta(\alpha+a+n)e^{2\pi inx}, ~x\in[0, 1].
\end{eqnarray*}
It follows from Propositions \ref{prop1.1} and \ref{toeplitz2} that $\Lambda_\alpha$ is a complete interpolation set for $V(\phi)$ if and only if $U_{11}$ is invertible on $\ell^2(\mathbb{N}^{-}).$ 
By \eqref{zak} and \eqref{fttheta}, we have
\begin{eqnarray}\label{symbol}
\Psi_{\alpha}(x)
&=&\mathcal{Z}\Theta(\alpha+a,-x)\nonumber\\
&=&e^{-2\pi i(\alpha+a) x}\mathcal{Z}\widehat{\Theta}(-x, -\alpha-a)\nonumber\\
&=&e^{-2\pi i(\alpha+a) x}\sum\limits_{n\in\mathbb{Z}}\widehat{\Theta}(n-x)e^{2\pi i(\alpha+a) n}\nonumber\\
&=&e^{-2\pi i(\alpha+a) x}\sum\limits_{n\in\mathbb{Z}}\widehat{\phi}(n-x)\Psi^{\dagger}(x-n)^{-1}e^{2\pi i(\alpha+a) n}\nonumber\\
&=&e^{-2\pi i(\alpha+a) x}\Psi^{\dagger}(x)^{-1}\sum\limits_{n\in\mathbb{Z}}\widehat{\phi}(n-x)e^{2\pi i(\alpha+a) n},~~a.e.,~ x\in[0, 1].
\end{eqnarray} 

We can not expect that $\Psi_{\alpha}$ is always continuous.
If $\Psi_{\alpha}$ is a piecewise continuous function with jumps, we will consider $\Psi_{\alpha}^{\#}$ the closed, continuous, and naturally oriented curve obtained from the range of $\Psi_\alpha.$
If $\Psi_\alpha$ is continuous, then 
$\Psi_{\alpha}^{\#}$ will be the same as $\Psi_\alpha$.
By Theorem \ref{toeplitz1}, we can conclude the following 
\begin{thm}\label{cis}
Let $\phi$ be a stable generator for $V(\phi)$. Suppose there exists $a\in [0,1)$ such that  $a+\mathbb{Z}$ is a complete interpolation set for $V(\phi)$ and $\Psi_{\alpha}$ is a piecewise continuous function on $[0,1]$. Then the set $\Lambda_{\alpha}=\{a+n:n\in\mathbb{N}_0\}\cup\{\alpha+a+n:n\in\mathbb{N}^{-}\}$ is a complete interpolation set for $V(\phi)$ if and only if $\Psi_{\alpha}^{\#}(x)\neq0$ for all $x\in[0, 1]$ and the winding number of $\Psi_{\alpha}^{\#}$ relative to zero is zero. 
\end{thm}
 
\section{Construction of complete interpolation sets for transversal-invariant spaces}

Throughout this paper,  $\lfloor x \rfloor$ denotes the floor function, $\lceil x \rceil$ denotes the ceiling function, and $\langle x\rangle$ denotes the fractional part of $x.$

Two measurable sets $A$ and $B$ in $\mathbb{R}$ are said to be translation congruent modulo $1$  (or $A$ is translation congruent to $B$ modulo $1$) if there exists a measurable bijection $\mu: A \to B$ such that
$$\mu(s)-s \in \mathbb{Z}, \text{ for each } s\in A.$$
Symbolically,  we write $A\equiv B \ (mod\ 1).$
Recall that a set $A\subseteq\mathbb{R}$  is said to be $1$-transversal if 
$\{A+n:n\in\mathbb{Z}\}$ forms a measurable partition of $\mathbb{R}$.
A  set $A$ is 1-transversal if and only if $A\equiv [0,1)\ (mod\ 1 )$ if and only if
\begin{equation}\label{eqn3.5}
\sum\limits_{n\in\mathbb{Z}}\chi_A(x+n)=1 ~\text{for a.e.}~ x\in\mathbb{R}.
\end{equation}
The condition \eqref{eqn3.5} implies that $\left\{\widehat{\chi_A}(\cdot-n):n\in\mathbb{Z}\right\}$
forms an orthonormal basis for $V(\widehat{\chi_A})$ and $\widehat{\chi_A}(n)=\delta_{n0}, n\in\mathbb{Z}$. It is easy to show that $\mathbb{Z}$ is a complete interpolation set for $V(\widehat{\chi_A})$.
Consequently, every $f\in V(\widehat{\chi_A})$ can be written as
$$f(x)=\sum\limits_{n\in\mathbb{Z}}f(n)\widehat{\chi_A}(x-n).$$
The sampling formula mentioned above is also proved in \cite{Torres}.

To determine the conditions on $\alpha$ for which the sample set $\mathbb{N}_0\cup\alpha+\mathbb{N}^{-}$ serves as a complete interpolation set for transversal invariant spaces, we need to prove the following lemmas which are essential for demonstrating the main result.

\begin{lem}\label{transversal}
Let $E\subseteq\mathbb{R},$ $F\subseteq[0, 1]$, and $\Delta=\big\{n\in\mathbb{Z}: n\in E+[0, 1]\big\}.$ If $E\equiv F\  (mod\ 1),$ then $$1-F=\bigcup\limits_{\lambda\in\Delta}\Big((\lambda-E)\cap [0, 1]\Big).$$
\end{lem}
\begin{pf}
Let $x\in\bigcup\limits_{\lambda\in\Delta}\Big((\lambda-E)\cap [0, 1]\Big).$ Then $x\in(\lambda-E)\cap [0, 1]$ for some $\lambda\in\Delta,$ which implies that $\lambda-x\in E$ and $1-x\in[0, 1].$ 
Since $E\equiv F\ (mod\ 1),$ there exists a unique $a\in F$ such that
\begin{eqnarray*}
\lambda-x\equiv a\ (mod\ 1).
\end{eqnarray*}
Hence $1-x\equiv a\ (mod\ 1)$ which implies that $x\in 1-F.$ 
Thus $\bigcup\limits_{\lambda\in\Delta}\Big((\lambda-E)\cap [0, 1]\Big)\subseteq1-F.$

Conversely, let $x\in1-F.$ Then $1-x\in F$
and hence there exists $y\in E$ such that $$y\equiv 1-x\ (mod\ 1)\implies y-1+x=n, \text{ for some } n\in\mathbb{Z}.$$ It is clear that $n+1=y+x\in\Delta$
and $x=n+1-y\in n+1-E.$ Since $x\in[0, 1],$ $x\in (n+1-E)\cap[0, 1].$
Thus $1-F\subseteq\bigcup\limits_{\lambda\in\Delta}\Big((\lambda-E)\cap[0, 1]\Big).$
\end{pf}
\begin{lem}\label{lineintegral}
Let $M, N\in\mathbb{Z}$, $A\in\mathbb{C}\setminus\{0\},$ $B\in\mathbb{R}$ and $\alpha\in\mathbb{R}-\left\{\frac{2l+1}{2(N-M)}: l\in\mathbb{Z}\right\}.$ 
\begin{itemize}
\item [$(i)$] For a curve $\gamma_1$ defined by $\gamma_1(t)=Ae^{2\pi iBt},~ a\leq t\leq b,$ we have
$$I_{\gamma_1}:=\dfrac{1}{2\pi i}\int\limits_{\gamma_1}\frac{dz}{z}=(b-a)B.$$
\item[$(ii)$] For a curve $\gamma_2$ defined by 
$\gamma_2(t)=A((1-t)e^{2\pi iM\alpha} + te^{2\pi iN\alpha}), ~0\leq t\leq1,$ we have 
$$I_{\gamma_2}:=\dfrac{1}{2\pi i}\int\limits_{\gamma}\frac{dz}{z}=(N-M)\alpha-\left\lfloor(N-M)\alpha+\frac{1}{2}\right\rfloor.$$ 
\end{itemize}
\end{lem}
\begin{pf}
The proof of $(i)$ is trivial.
 To prove $(ii)$, let $N=M+n,$ for some $n\in\mathbb{Z}.$ Then we rewrite $\gamma_2$ as
\begin{eqnarray*}
\gamma_2(t)&=&A\left\{(1-t)e^{2\pi iM\alpha} + te^{2\pi i(M+n)\alpha}\right\}\\
&=&Ae^{2\pi i(M+\frac{n}{2})\alpha} \left(e^{-\pi in\alpha}+t(e^{\pi in\alpha}-e^{-\pi in\alpha})\right)\\
&=&Ae^{2\pi i(M+\frac{n}{2})\alpha}\left(\cos{\pi n\alpha}-i(1-2t)\sin{\pi n\alpha}\right)\\
&=&Ae^{2\pi i(M+\frac{n}{2})\alpha}\cos{\pi n\alpha}[1-i(1-2t)\tan{\pi n\alpha}]\\
&=&Ae^{2\pi i(M+\frac{n}{2})\alpha}\cos{\pi n\alpha}\dfrac{[1+(1-2t)^2\tan^2{\pi n\alpha}]}{1+i(1-2t)\tan{\pi n\alpha}}.
\end{eqnarray*}
 Differentiating $\gamma_2$ with respect to $t,$ we get
 $\gamma_2^{\prime}(t)=2iAe^{2\pi i(M+\frac{n}{2})\alpha}\sin{\pi n\alpha}.$
 Therefore, 
 \begin{eqnarray*}
I_{\gamma_2}&=&\dfrac{\tan{\pi n\alpha}}{\pi}\int\limits_{0}^{1}\dfrac{1+i(1-2t)\tan{\pi n\alpha}}{1+(1-2t)^2\tan^2{\pi n\alpha}}~dt\\
&=&\dfrac{\tan{\pi n\alpha}}{\pi}\left\{\int\limits_{0}^{1}\dfrac{dt}{1+(1-2t)^2\tan^2{\pi n\alpha}}+i\tan{\pi n\alpha}\int\limits_{0}^{1}\dfrac{(1-2t)dt}{1+(1-2t)^2\tan^2{\pi n\alpha}}\right\}\\
&=&\dfrac{\tan{\pi n\alpha}}{\pi}\left\{\frac{1}{2}\int\limits_{-1}^{1}\dfrac{du}{1+u^2\tan^2{\pi n\alpha}}+\frac{i\tan{\pi n\alpha}}{2}\int\limits_{-1}^{1}\dfrac{u}{1+u^2\tan^2{\pi n\alpha}}du\right\}\\
&=&\dfrac{\tan^{-1}(\tan{\pi n\alpha})}{\pi}=\dfrac{\tan^{-1}({\tan{\pi (N-M)\alpha}})}{\pi}\\
&=&(N-M)\alpha-\left\lfloor(N-M)\alpha+\frac{1}{2}\right\rfloor.
\end{eqnarray*}
\end{pf}
\begin{lem}\label{flr}
Let $l\in\mathbb{N}$ and $r\in\mathbb{Z}\setminus\{0\}$ such that $\tfrac{l}{r}\in\mathbb{Z}.$ If $\alpha\in\left[\frac{n-1}{2l}, \frac{n}{2l}\right),$ $n\in\mathbb{Z},$ then 
\begin{eqnarray}\label{flr1}
\left\lfloor\alpha r+\frac{1}{2}\right\rfloor=\begin{cases}
\left\lfloor\frac{(n-1)r}{2l}+\frac{1}{2}\right\rfloor  &\text{ if } r>0,\\\\
\left\lfloor\frac{nr}{2l}+\frac{1}{2}\right\rfloor  &\text{ if } r<0.
\end{cases}    
\end{eqnarray}
\end{lem}
\begin{pf}
If $r>0$ and $\frac{n-1}{2l}\leq\alpha <\frac{n}{2l}$, then
\begin{eqnarray*}
\frac{(n-1)r}{2l}+\frac{1}{2}&\leq&\alpha r+\frac{1}{2}<\frac{nr}{2l}+\frac{1}{2},
\end{eqnarray*}
which implies that
\begin{eqnarray*}
\frac{nr+l}{2l}-\frac{r}{2l}&\leq&\alpha r+\frac{1}{2}<\frac{nr+l}{2l}.
\end{eqnarray*}
It is clear that when $\dfrac{nr+l}{2l}\in\mathbb{Z},$ then \eqref{flr1} is always true.
When $\dfrac{nr+l}{2l}\not\in\mathbb{Z},$
then there does not exists any integer $m$ such that $\frac{nr+l}{2l}-\frac{r}{2l}< m<\frac{nr+l}{2l}.$
(If it exists, then $(2m-1)\dfrac{l}{r}<n<(2m-1)\dfrac{l}{r}+1$ which is impossible).
Hence \eqref{flr1} is proved for $r>0$.
Similarly, we can prove the result  for $r<0$.
\end{pf}

Let $E$ be a 1-transversal set satisfying the following properties:
\begin{itemize}
\item [$(P_1)$] $E=\bigcup\limits_{k=1}^{L}E_k$, $E_k$'s are disjoint intervals,
\item [$(P_2)$]  $E_k\equiv F_k~ (mod~ 1),$ where
$F_k=[a_{k}, a_{k+1}]$ with $0=a_1<a_2<\dots<a_{L+1}=1.$ 
\end{itemize}

Let us define
$$\Delta_k:=\Big\{\lambda\in\mathbb{Z}: \lambda\in E_k+[0, 1]\Big\}\text{ and } \mu_k=\#(\Delta_k),~1\leq k\leq L.$$ It is clear that $1\leq\mu_k\leq2$. We assume that $\Delta_{k}=\{\lambda_{k1}, \lambda_{k\mu_k}:\lambda_{k1}\leq\lambda_{k\mu_k}\},$ for all $1\leq k\leq L.$ If $\mu_k=2,$ then $\lambda_{k2}=1+\lambda_{k1}.$

Let $\nu=\lcm\{|\lambda_{k,\mu_k}-\lambda_{k1}|, ~|\lambda_{k-1,1}-\lambda_{k\mu_k}|\in\mathbb{N}: 1\leq k\leq L\}.$  Here $\lambda_{01}=\lambda_{L1}+1$.
Let $\Omega=\Big\{ k\in\{1,2,\dots L\}:\mu_k=2\Big\}$ and $\#(\Omega)=\rho.$ 
Let us define 
$$g(k)=\begin{cases}
1&\text{ if } \lambda_{k-1,1}-\lambda_{k1}>0,\\
0&\text{ if } \lambda_{k-1,1}-\lambda_{k1}<0,
\end{cases}$$ 
\begin{eqnarray*}
f(n)&=&-\rho\left\lfloor\frac{n-1}{2\nu}+\frac{1}{2}\right\rfloor-\sum\limits_{k\in\Omega}\left\lfloor\frac{n-g(k)}{2\nu}(\lambda_{k-1,1}-\lambda_{k1}-1)+\frac{1}{2}\right\rfloor\\
&&-\sum\limits_{k\in\Omega^\complement}\left\lfloor\frac{n-g(k)}{2\nu}(\lambda_{k-1,1}-\lambda_{k1})+\frac{1}{2}\right\rfloor,  
\end{eqnarray*}
and
$$
A=\Bigg\{n\in\left\{-2v\left\lceil \frac{L+\rho}{2}\right\rceil,\dots, 2v\left\lceil \frac{L+\rho}{2}\right\rceil\right\} :
f(n)=0\Bigg\}.   
$$
\begin{thm}\label{general}
Let $E$ be a 1-transversal set satisfying the conditions $P_1$ and $P_2$. Let $\Delta_k$ and $A$ be defined as above. Then $\mathbb{N}_0\cup\alpha+\mathbb{N}^{-}$ is a complete interpolation set for $V(\widehat{\chi_E})$ if and only if
\begin{itemize}
    \item [$(i)$] $|\alpha|\leq\frac{L+\rho}{2}$,
    \item[$(ii)$] $\alpha\in\bigcup\limits_{n\in A}\left[\dfrac{n-1}{2\nu}, \dfrac{n}{2\nu}\right)\setminus{G},$ where 
$$G=\left\{\dfrac{2l+1}{2(\lambda_{k\mu_k}-\lambda_{k1})}, \dfrac{2l+1}{2(\lambda_{k-1,1}-\lambda_{k\mu_k})}\in\mathbb{R}: l\in\mathbb{Z}, ~1\leq k\leq L\right\}.$$
\end{itemize}
\end{thm}
\begin{pf}
For $\phi=\widehat{\chi_{E}},$ we have $\Psi^{\dagger}(x)=1.$ Since $E_k$'s are  disjoint sets,
it follows from \eqref{symbol} that
\begin{eqnarray}\label{transversalpsi}
\Psi_{\alpha}(x)
&=&e^{-2\pi i\alpha x}\sum\limits_{k=1}^{L}\sum\limits_{n\in\mathbb{Z}}\chi_{E_k}(n-x)e^{2\pi i\alpha n}\nonumber\\
&=&e^{-2\pi i\alpha x}\sum\limits_{k=1}^{L}\sum\limits_{n\in\Delta_k}\chi_{E_k}(n-x)e^{2\pi i\alpha n}\nonumber\\
&=&e^{-2\pi i\alpha x}\sum\limits_{k=1}^{L}\sum\limits_{i=1}^{\mu_k}\chi_{E_k}(\lambda_{ki}-x)e^{2\pi i\alpha\lambda_{ki}},~a.e.,~x\in[0, 1].
\end{eqnarray}
Since $E$ is a 1-transversal set, \eqref{transversalpsi} becomes
\begin{eqnarray*}
\Psi_{\alpha}(x)=e^{-2\pi i\alpha x}\begin{cases}
e^{2\pi i\alpha\lambda_{11}} &\text{ if } x\in(\lambda_{11}-E_1^\circ)\cap(0, 1),\\
e^{2\pi i\alpha\lambda_{1\mu_1}} &\text{ if } x\in(\lambda_{1\mu_1}-E_1^\circ)\cap(0, 1),\\
\hspace{.5cm}\vdots\\
e^{2\pi i\alpha\lambda_{L1}} &\text{ if } x\in(\lambda_{L1}-E_L^\circ)\cap(0, 1),\\
e^{2\pi i\alpha\lambda_{L\mu_L}} &\text{ if } x\in(\lambda_{L\mu_L}-E_L^\circ)\cap(0, 1).
\end{cases}
\end{eqnarray*}
Here $E^\circ$ denotes the interior of the set $E$ and the equality holds pointwise.
Notice that if $\lambda_{k1}\neq\lambda_{k\mu_k}$, then $\Big((\lambda_{k1}-E_k^\circ)\cap(0,1)\Big)\bigcap\Big((\lambda_{k\mu_k}-E_k^\circ)\cap(0,1)\Big)=\emptyset.$
By Lemma \ref{transversal}, we have  $$\bigcup\limits_{i=1}^{\mu_k}\Big((\lambda_{ki}-E_k^\circ)\cap(0, 1)\Big)=1-F_k^\circ=(1-a_{k+1}, 1-a_{k}),~1\leq k\leq L.$$
Therefore, if $x\in(\lambda_{k1}-E_k^\circ)\cap(0,1),$ then there exists $s_k\in1-F_k$ such that either $x\in(1-a_{k+1}, s_k]$ or $x\in(s_k, 1-a_k).$
Hence we can write
\begin{eqnarray*}
\Psi_{\alpha}(x)=e^{-2\pi i\alpha x}\begin{cases}
e^{2\pi i\alpha\lambda_{1\mu_1}} &\text{ if } x\in(s_1, 1-a_1),\\
e^{2\pi i\alpha\lambda_{11}} &\text{ if } x\in(1-a_2, s_1],\\
\hspace{0.5cm}\vdots\\
e^{2\pi i\alpha\lambda_{L\mu_L}} &\text{ if } x\in(s_L, 1-a_L),\\
e^{2\pi i\alpha\lambda_{L1}} &\text{ if } x\in (1-a_{L+1}, s_L].
\end{cases}
\end{eqnarray*}
Since $a_1=0$ and $a_{L+1}=1,$ we rewrite $\Psi_{\alpha}$ as
\begin{eqnarray*}
\Psi_{\alpha}(x)=e^{-2\pi i\alpha x}\begin{cases}
e^{2\pi i\alpha\lambda_{L1}} &\text{ if } x\in(0, s_L]\\
e^{2\pi i\alpha\lambda_{L\mu_L}} &\text{ if } x\in(s_L, 1-a_L),\\
e^{2\pi i\alpha\lambda_{L-1,1}} &\text{ if } x\in(1-a_{L}, s_{L-1}]\\
e^{2\pi i\alpha\lambda_{L-1,\mu_{L-1}}} &\text{ if } x\in(s_{L-1}, 1-a_{L-1}),\\
\hspace{0.5cm}\vdots\\
e^{2\pi i\alpha\lambda_{21}} &\text{ if } x\in(1-a_3, s_2],\\
e^{2\pi i\alpha\lambda_{2\mu_2}} &\text{ if } x\in(s_2, 1-a_2),\\
e^{2\pi i\alpha\lambda_{11}} &\text{ if } x\in(1-a_2, s_1],\\
e^{2\pi i\alpha\lambda_{1\mu_1}} &\text{ if } x\in(s_1, 1).
\end{cases}
\end{eqnarray*}
We do not always expect the function $\Psi_\alpha(x)$ to be continuous on $[0, 1]$. The function $\Psi_\alpha(x)$ may have jump discontinuities at $x=s_k,~1-a_k,$ for $k=L, L-1,\dots,1.$ To construct a closed curve  $\Psi_\alpha^{\#},$ we define the line segment joining $e^{-2\pi i\alpha s_k}e^{2\pi i\alpha\lambda_{k1}}$ to $e^{-2\pi i\alpha s_k}e^{2\pi i\alpha\lambda_{k\mu_k}}$ as
$$S_k(t)=e^{-2\pi i\alpha s_k}\Big((1-t)e^{2\pi i\alpha\lambda_{k1}}+te^{2\pi i\alpha\lambda_{k\mu_k}}\Big),~t\in[0, 1],$$
and the line segment joining $e^{-2\pi i\alpha (1-a_k)}e^{2\pi i\alpha\lambda_{k\mu_k}}$ to $e^{-2\pi i\alpha (1-a_k)}e^{2\pi i\alpha\lambda_{k-1,1}}$  as
$$\widetilde{S_k}(t)=e^{-2\pi i\alpha(1-a_k)}\Big((1-t)e^{2\pi i\alpha\lambda_{k\mu_k}}+te^{2\pi i\alpha\lambda_{k-1,1}}\Big),~t\in[0, 1],$$
for each $k=L, L-1,\dots,1.$
It is easy to check that 
\begin{itemize}
    \item[$(i)$] $S_k(t)=0$ if and only if $t=\frac{1}{2},$ $\alpha\in\left\{\frac{2l+1}{2(\lambda_{k\mu_k}-\lambda_{k1})}\in\mathbb{R}: l\in\mathbb{Z}\right\}$;\\
    \item[$(ii)$] $\widetilde{S_k}(t)=0$ if and only if $t=\frac{1}{2},$ $\alpha\in\left\{\frac{2l+1}{2(\lambda_{k-1,1}-\lambda_{k\mu_k})}\in\mathbb{R}: l\in\mathbb{Z}\right\}$.
\end{itemize}
Further, we define the curves
$$C_k(x)=e^{-2\pi i\alpha x}e^{2\pi i\alpha\lambda_{k1}},~x\in[1-a_{k+1}, s_k],$$ 
and  
$$\widetilde{C_k}(x)=e^{-2\pi i\alpha x}e^{2\pi i\alpha\lambda_{k\mu_k}},~x\in[s_k, 1-a_k],$$
for each $k=L, L-1,\dots,1.$
Now the curve $\Psi_\alpha^{\#}$ defined by
\begin{eqnarray*}
\Psi_\alpha^{\#}:=C_L+S_L+ \widetilde{C_L}+\widetilde{S_L}+\dots+C_1+S_1+\widetilde{C_1}+\widetilde{S_1}
\end{eqnarray*}
is closed. It is clear that
the curve $\Psi_\alpha^{\#}$ does not pass through the origin
if and only if
$\alpha\in\mathbb{R}-G.$
By Lemma \ref{lineintegral}, we have
\begin{align*}
I_{C_k}&=-\alpha(s_k-1+a_{k+1}),~
I_{S_k}=\alpha(\lambda_{k\mu_k}-\lambda_{k1})-\left\lfloor\alpha(\lambda_{k\mu_k}-\lambda_{k1})+\frac{1}{2}\right\rfloor,\\
I_{\widetilde{C_{k}}}&=-\alpha(1-a_k-s_k),~
I_{\widetilde{S_k}}=\alpha(\lambda_{k-1,1}-\lambda_{k\mu_k})-\left\lfloor\alpha(\lambda_{k-1,1}-\lambda_{k\mu_k})+\frac{1}{2}\right\rfloor,
\end{align*}
for $k=L, L-1,\dots,1.$
Hence the winding number of $\Psi_\alpha^{\#}$ relative to zero is 
\begin{align}\label{floor}
\textrm{Ind}(\Psi_\alpha^{\#},0)
&=\sum\limits_{k=L}^{1}\left(I_{C_k}+I_{\widetilde{C_{k}}}+I_{S_k}+I_{\widetilde{S_k}}\right)\nonumber\\
&=\alpha\sum\limits_{k=L}^{1}\left(-s_k+1-a_{k+1}-1+a_k+s_k+\lambda_{k\mu_k}-\lambda_{k1}+\lambda_{k-1,1}-\lambda_{k\mu_k}\right)\nonumber\\
&\hspace{2cm}-\sum\limits_{k=L}^{1}\left(\left\lfloor\alpha(\lambda_{k\mu_k}-\lambda_{k1})+\dfrac{1}{2}\right\rfloor+\left\lfloor\alpha(\lambda_{k-1,1}-\lambda_{k\mu_k})+\dfrac{1}{2}\right\rfloor\right)\nonumber\\
&=-\sum\limits_{k=L}^{1}\Big(\left\lfloor\alpha(\lambda_{k\mu_k}-\lambda_{k1})+\frac{1}{2}\right\rfloor+\left\lfloor\alpha(\lambda_{k-1,1}-\lambda_{k\mu_k})+\frac{1}{2}\right\rfloor\Big)\\
&=-\sum\limits_{k\in\Omega}\left\lfloor\alpha(\lambda_{k2}-\lambda_{k1})+\frac{1}{2}\right\rfloor-\sum\limits_{k\in\Omega}\left\lfloor\alpha(\lambda_{k-1,1}-\lambda_{k2})+\frac{1}{2}\right\rfloor\nonumber\\
&\hspace{2cm}-\sum\limits_{k\in\Omega^\complement}\left\lfloor\alpha(\lambda_{k-1,1}-\lambda_{k1})+\frac{1}{2}\right\rfloor\nonumber\\
&=-\rho\left\lfloor\alpha+\frac{1}{2}\right\rfloor-\sum\limits_{k\in\Omega^\complement}\left\lfloor\alpha(\lambda_{k-1,1}-\lambda_{k1})+\frac{1}{2}\right\rfloor\nonumber\\
&\hspace{5cm}-\sum\limits_{k\in\Omega}\left\lfloor\alpha(\lambda_{k-1,1}-\lambda_{k1}-1)+\frac{1}{2}\right\rfloor.\label{flr2}
\end{align}
It follows from \eqref{floor} that
\begin{align}
\textrm{Ind}(\Psi_\alpha^{\#},0)
&=-\sum\limits_{k=L}^{1}\left\{\alpha(\lambda_{k\mu_k}-\lambda_{k1})+\frac{1}{2}+\alpha(\lambda_{k-1,1}-\lambda_{k\mu_k})+\frac{1}{2}\right\}\nonumber\\
&\hspace{1.5cm}+\sum\limits_{k=L}^{1}\left\{\left\langle\alpha(\lambda_{k\mu_k}-\lambda_{k1})+\frac{1}{2}\right\rangle+\left\langle\alpha(\lambda_{k-1,1}-\lambda_{k\mu_k})+\frac{1}{2}\right\rangle\right\}\nonumber\\
&=-\alpha-L+\sum\limits_{k=L}^{1}\left\{\left\langle\alpha(\lambda_{k\mu_k}-\lambda_{k1})+\frac{1}{2}\right\rangle+\left\langle\alpha(\lambda_{k-1,1}-\lambda_{k\mu_k})+\frac{1}{2}\right\rangle\right\}\nonumber
\end{align}
\begin{align}
&=-\alpha-L+\frac{L-\rho}{2}+\rho\left\langle\alpha+\frac{1}{2}\right\rangle\nonumber\\
&\hspace{1.5cm}+\sum\limits_{k\in\Omega}\left\langle\alpha(\lambda_{k-1,1}-\lambda_{k1}-1)+\frac{1}{2}\right\rangle+\sum\limits_{k\in\Omega^\complement}\left\langle\alpha(\lambda_{k-1,1}-\lambda_{k1})+\frac{1}{2}\right\rangle
\nonumber\\
&=-\alpha-\frac{L+\rho}{2}+\rho\left\langle\alpha+\frac{1}{2}\right\rangle\nonumber\\
&\hspace{1.5cm}+\sum\limits_{k\in\Omega}\left\langle\alpha(\lambda_{k-1,1}-\lambda_{k1}-1)+\frac{1}{2}\right\rangle+\sum\limits_{k\in\Omega^\complement}\left\langle\alpha(\lambda_{k-1,1}-\lambda_{k1})+\frac{1}{2}\right\rangle.\label{fractional}
\end{align} 
By using the fact $0\leq\langle x\rangle<1,$ we can conclude from \eqref{fractional} that $$\textrm{Ind}(\Psi_\alpha^{\#},0)\neq0,\text{ if }|\alpha|>\frac{L+\rho}{2}.$$
Suppose $\alpha\in\left[\dfrac{n-1}{2\nu}, \dfrac{n}{2\nu}\right)\setminus{G},$ where $n\in\mathbb{Z}$ such that $|n|\leq2\nu\left\lceil\dfrac{L+\rho}{2} \right\rceil.$ By Lemma \ref{flr}, we obtain from \eqref{flr2} that
\begin{eqnarray}
\textrm{Ind}(\Psi_\alpha^{\#},0)
=f(n)=0 \text{ if and only if } n\in A.
\end{eqnarray}
Our result now follows from Theorem \ref{cis}.
\end{pf}
\begin{cor}\label{rbthm}
Let $E$ be a 1-transversal set satisfying the conditions $P_1$ and $P_2$ and $\Delta_k$ defined as above. Then $\left\{e^{2\pi i\xi\cdot}:\xi\in\mathbb{N}_0\cup\alpha+\mathbb{N}^{-}\right\}$ is a Riesz basis for $L^2(E)$ if and only if $\alpha$ satisfies $(i)$ and $(ii)$ of Theorem \ref{general}.
\end{cor}
\begin{pf}
From \eqref{zak}, we have
\begin{eqnarray*}
\mathcal{Z}\phi(\xi,y)=e^{2\pi i\xi y}\sum\limits_{k\in\mathbb{Z}}\chi_E(y-k)e^{-2\pi ik\xi}
=e^{2\pi i \xi y} \text{ if } \xi\in \mathbb{N}_0\cup\alpha+\mathbb{N}^{-}\text{ and } y\in E.
\end{eqnarray*}
Now our result follows from Proposition \ref{zakthm} and Theorem \ref{general}.
\end{pf}
\begin{cor}
Let $E$ be a 1-transversal set satisfying the conditions $P_1$ and $P_2$ and $\Delta_k$ defined as above. 
\begin{itemize}
    \item [$(i)$] If $\alpha\in\left[-\frac{1}{2\nu}, \frac{1}{2v}\right)\setminus{G}$, then $\Lambda_\alpha$ is  a complete interpolation set for $V(\widehat{\chi_E})$.
    \item[$(ii)$] If $\alpha\in\left[m-\frac{1}{2\nu}, m+\frac{1}{2v}\right),$ $m\in\mathbb{Z}-\{0\}$, then $\Lambda_\alpha$ is not a complete interpolation set for $V(\widehat{\chi_E})$.
\end{itemize}
\end{cor}
\begin{pf}
If $n=2m\nu,$ $m\in\mathbb{Z},$ then
\begin{align*}
f(n)&=-\rho\left\lfloor m-\frac{1}{2\nu}+\frac{1}{2}\right\rfloor-\sum\limits_{k\in\Omega^\complement}\left\lfloor m(\lambda_{k-1,1}-\lambda_{k1})-\frac{g(k)}{2\nu}(\lambda_{k-1,1}-\lambda_{k1})+\frac{1}{2}\right\rfloor\\
&-\sum\limits_{k\in\Omega}\left\lfloor m(\lambda_{k-1,1}-\lambda_{k1}-1)-\frac{g(k)}{2\nu}(\lambda_{k-1,1}-\lambda_{k1}-1)+\frac{1}{2}\right\rfloor\\
&=-m\rho-\rho\left\lfloor -\frac{1}{2\nu}+\frac{1}{2}\right\rfloor-m\sum\limits_{k\in\Omega}(\lambda_{k-1,1}-\lambda_{k1}-1)-m\sum\limits_{k\in\Omega^\complement}(\lambda_{k-1,1}-\lambda_{k1})\\
&-\sum\limits_{k\in\Omega}\left\lfloor -\frac{g(k)}{2\nu}(\lambda_{k-1,1}-\lambda_{k1}-1)+\frac{1}{2}\right\rfloor-\sum\limits_{k\in\Omega^\complement}\left\lfloor -\frac{g(k)}{2\nu}(\lambda_{k-1,1}-\lambda_{k1})+\frac{1}{2}\right\rfloor\\
&=-m-\sum\limits_{k\in\Omega}\left\lfloor -\frac{g(k)}{2\nu}(\lambda_{k-1,1}-\lambda_{k1}-1)+\frac{1}{2}\right\rfloor-\sum\limits_{k\in\Omega^\complement}\left\lfloor -\frac{g(k)}{2\nu}(\lambda_{k-1,1}-\lambda_{k1})+\frac{1}{2}\right\rfloor\\
&=-m,
\end{align*}
because $\frac{1}{2}-\frac{1}{2\nu},~-\frac{g(k)}{2\nu}(\lambda_{k-1,1}-\lambda_{k1}-1)+\frac{1}{2},~-\frac{g(k)}{2\nu}(\lambda_{k-1,1}-\lambda_{k1})+\frac{1}{2}\in[0, 1).$
Similarly,  we can show that if $n=2m\nu+1,$ $m\in\mathbb{Z},$ then $f(n)=-m$.
Hence our result follows from Theorem $\ref{general}.$
\end{pf}
 
We will illustrate our theorem by examining specific transversal sets. We refer to \cite{Dai, Shukla} for examples of transversal sets.
\begin{example}
We shall carry out the computation for
$E=\left[-\frac{1}{2}, \frac{1}{2}\right).$ In this case, we have 
$$L=1,~\Delta_1=\{\lambda_{11}=0, \lambda_{12}=1\},~\rho=1, \text{ and } \nu=1.$$
Now it is easy to conclude from Theorem \ref{general} and Corollary \ref{rbthm} that
 $\mathbb{N}_0\cup\alpha+\mathbb{N}^{-}$ is a complete interpolation set for $V(\sinc)$ if and only if  $\left\{e^{2\pi i\xi\cdot}:\xi\in\mathbb{N}_0\cup\alpha+\mathbb{N}^{-}\right\}$ is a Riesz basis for $L^2[-\frac{1}{2}, \frac{1}{2})$ if and only if $|\alpha|<\frac{1}{2}.$
\end{example}
\begin{example}
Consider the Littlewood-Paley wavelet
\begin{eqnarray*}
\phi(x):=\dfrac{\sin\pi x}{\pi x}(2\cos\pi x-1).
\end{eqnarray*}
Its Fourier transform is given by
\begin{equation*}
 \widehat{\phi}(w)= \left\{\begin{array}{cc}
1& \mbox{if $w\in E$,}\\
0 &\mbox{otherwise.}
\end{array} \right.
\end{equation*}
where $E=[-1,-1/2)\cup(1/2,1]$. Since $E$ is a 1-transversal set,  every $f\in V(\phi)$ can be written as
$$f(x)=\sum\limits_{n\in\mathbb{Z}}f(n)\dfrac{\sin\pi(x-n)}{\pi(x-n)}[2\cos \pi(x-n)-1].$$
It is clear that
 $\left[-1, -\frac{1}{2}\right)\equiv\left[0,\frac{1}{2}\right)\ (mod\ 1)$ and $\left(\frac{1}{2}, 1\right]\equiv\left(\frac{1}{2}, 1\right]\ (mod\ 1)$.
In this case,  we have
$$L=2,~\lambda_{11}=0,~\lambda_{21}=1,~\lambda_{01}=\lambda_{21}+1=2,~\rho=0,$$ $$\nu=\lcm\{|\lambda_{01}-\lambda_{11}|=2, |\lambda_{11}-\lambda_{21}|=1\}=2,\text{ and }G=\left\{\frac{2l+1}{4}, \frac{2l+1}{2}:l\in\mathbb{Z}\right\}.$$ Further,
$$
f(n)
=-\left\lfloor\frac{n}{2}\right\rfloor-\left\lfloor-\frac{n}{4}+\frac{1}{2}\right\rfloor=0 \text{ if and only if } n\in A=\{-2, 0, 1, 3\}.$$
By Theorem \ref{general}, the sample set $\mathbb{N}_0\cup\alpha+\mathbb{N}^{-}$ is a complete interpolation set for $V(\phi)$ if and only if $\alpha\in\left(-\frac{1}{4}, \frac{1}{4}\right)\cup\left(-\frac{3}{4}, -\frac{2}{4}\right)\cup\left(\frac{2}{4}, \frac{3}{4}\right).$
\end{example}
\begin{example}
Consider the Journe transversal set $$E=\left[2, \frac{16}{7}\right)\cup\left[\frac{2}{7}, \frac{1}{2}\right)\cup\left[-\frac{1}{2}, -\frac{2}{7}\right)\cup\left[-\frac{16}{7}, -2\right).$$
It is clear that $\left[2, \frac{16}{7}\right)\equiv\left[0, \frac{2}{7}\right)\ (mod\ 1),$  $\left[-\frac{1}{2}, -\frac{2}{7}\right)\equiv\left[\frac{1}{2}, \frac{5}{7}\right)\ (mod\ 1),$ and $\left[-\frac{16}{7}, -2\right)\equiv\left[\frac{5}{7}, 1\right)\ (mod\ 1).$
In this case, we have $$L=4,~\lambda_{11}=3,~\lambda_{21}=2,~\lambda_{31}=0,~\lambda_{41}=-2,~\lambda_{01}=\lambda_{41}+1=-1,~\rho=0,$$ 
 $$~\nu=4,\text{ and }G=\left\{\frac{2l+1}{8}, \frac{2l+1}{2}, \frac{2l+1}{4}:l\in\mathbb{Z}\right\}.$$ 
Further,
$$
f(n)
=-\left\lfloor-\frac{n}{2}+\frac{1}{2}\right\rfloor-\left\lfloor\frac{n}{8}+\frac{3}{8}\right\rfloor-2\left\lfloor\frac{n}{4}+\frac{1}{4}\right\rfloor=0$$ 
if and only if  $n\in A=\{-9, -5, -3, 0, 1, 4, 6, 10\}.$
By Theorem \ref{general}, $\mathbb{N}_0\cup\alpha+\mathbb{N}^{-}$ is a complete interpolation set for $V(\widehat{\chi_E})$ if and only if $\alpha\in\left(-\frac{10}{8}, -\frac{9}{8}\right)\cup\left(-\frac{6}{8}, -\frac{5}{8}\right)\cup\left(-\frac{4}{8}, -\frac{3}{8}\right)\cup\left(-\frac{1}{8}, \frac{1}{8}\right)\cup\left(\frac{3}{8}, \frac{4}{8}\right)\cup\left(\frac{5}{8}, \frac{6}{8}\right)\cup\left(\frac{9}{8}, \frac{10}{8}\right).$ 
\end{example}
\begin{example}
Let us choose a non-symmetric transversal set
$$E=\left[3, \dfrac{16}{5}\right)\cup\left[\dfrac{1}{5}, \dfrac{1}{3}\right)\cup\left[\dfrac{4}{3}, \dfrac{3}{2}\right)\cup\left[-\dfrac{1}{2}, -\dfrac{1}{3}\right)\cup\left[-\dfrac{4}{3}, -1\right).$$
Then  we have
$$\left[3, \dfrac{16}{5}\right)\equiv\left[0, \dfrac{1}{5}\right)\ (mod\ 1),~\left[\dfrac{4}{3}, \dfrac{3}{2}\right)\equiv\left[\dfrac{1}{3}, \dfrac{1}{2}\right)\ (mod\ 1),$$
$$\left[-\dfrac{1}{2}, -\dfrac{1}{3}\right)\equiv\left[\dfrac{1}{2}, \dfrac{2}{3}\right)\ (mod\ 1), \text{ and }\left[-\dfrac{4}{3}, -1\right)\equiv\left[\dfrac{2}{3}, 1\right)\ (mod\ 1).$$
In this case, $L=5,$ $\rho=0,$ $\lambda_{11}=4,$ $\lambda_{21}=1,$ $\lambda_{31}=2,$ $\lambda_{41}=0,$ $\lambda_{51}=-1,$$\lambda_{01}=\lambda_{51}+1=0,$ $\nu=12$ and $G=\left\{\dfrac{2l+1}{8}, \dfrac{2l+1}{6}, \dfrac{2l+1}{2}, \dfrac{2l+1}{4}: l\in\mathbb{Z}\right\}.$
Further, $$f(n)=-\left\lfloor-\dfrac{n}{6}+\dfrac{1}{2}\right\rfloor-\left\lfloor\dfrac{n}{8}+\dfrac{3}{8}\right\rfloor-\left\lfloor-\dfrac{n}{24}+\dfrac{1}{2}\right\rfloor-\left\lfloor\dfrac{n}{12}+\dfrac{5}{12}\right\rfloor-\left\lfloor\dfrac{n}{24}+\dfrac{11}{24}\right\rfloor=0$$ 
if and only if $A=\{0, \pm1, \pm2, 3, -4, \pm5, 6, -9, \pm10, \pm11, 12, -15, \pm16, \pm17, 18, -27, 28 \}.$

By Theorem \ref{general}, $\mathbb{N}_0\cup\alpha+\mathbb{N}^{-}$ is a complete interpolation set for $V(\widehat{\chi_E})$ if and only if $\alpha\in\left(-\frac{1}{8}, \frac{1}{8}\right)\cup\left(-\frac{1}{4}, -\frac{1}{6}\right)\cup\left(\frac{1}{6}, \frac{1}{4}\right)\cup\left(-\frac{1}{2}, -\frac{3}{8}\right)\cup\left(\frac{3}{8}, \frac{1}{2}\right)\cup\left(-\frac{3}{4}, -\frac{5}{8}\right)\cup\left(\frac{5}{8}, \frac{3}{4}\right)\cup\left(-\frac{7}{6}, -\frac{9}{8}\right)\cup\left(\frac{9}{8}, \frac{7}{6}\right).$ 
\end{example}

\section{Some new properties of exponential splines}
The $B$-splines $Q_m$ are defined inductively as follows:
\begin{align*}
Q_1(x)=\chi_{[0,1)}(x) \text{ and }Q_{m+1}(x)=\int\limits_{-\infty}^{\infty}Q_m(x-t)Q_1(t)dt,
~m\geq1.
\end{align*}
The $B$-spline $Q_m$ is compactly supported on the interval $\left[0,m\right]$
and its Fourier transform is given by $\widehat{Q_{m}}(w)=e^{-\pi imw}\left(\sinc{ w}\right)^m,$ where $\sinc{y}=\dfrac{\sin{\pi y}}{\pi y}.$
The collection $\big\{Q_m(\cdot-n):n\in\mathbb{Z}\big\}$ forms a Riesz basis for  the shift-invariant spline space $V(Q_m)$. 

The exponential spline $\Phi_{m-1}(\beta,t)$ of degree $m-1$ to the base $t$ is defined by
\begin{equation}
\Phi_{m-1}(\beta,t)=\sum\limits_{k\in\mathbb{Z}}t^kQ_{m}(\beta-k),~ 0\leq \beta<1,~t\neq0,1.
\end{equation}
Schoenberg \cite{Schoenberg} provided a following recurrence relation for exponential splines: 
\begin{equation*}
\int\limits_{\beta}^{\beta+1}\Phi_{m-2}(u,t)~du=t\Phi_{m-1}(\beta,t).
\end{equation*}

To analyze the zeros of exponential splines with respect to the variable $t$, we need to introduce a new recurrence relation. The above recurrence relation is inadequate for characterizing a complete interpolation set for shift-invariant spline spaces. Thus, we will develop and use our new recurrence relation to better understand and analyze the behavior of these splines.

The Bernoulli numbers $B_n, n\geq 0$ and the Euler polynomials $E_n(x), n\geq 0$ are  respectively defined by
\begin{eqnarray}
B_n=-\frac{1}{n+1}\sum\limits_{k=0}^{n-1}\binom{n+1}{k}B_k, ~B_0=1,\nonumber
\end{eqnarray}
and
\begin{equation}
E_n(x)=\frac{1}{n+1}\sum\limits_{k=1}^{n+1}\binom{n+1}{k}(2-2^{k+1})B_kx^{n+1-k}.\nonumber
\end{equation}
For further details, we refer to \cite{Aci}.
The first few Bernoulli numbers and Euler polynomials are
$$B_0=1,~B_1=-\frac{1}{2},~B_2=\frac{1}{6},~B_3=0,~B_4=-\frac{1}{30},~B_5=0,~B_6=\frac{1}{42},$$
and
$$E_0(x)=1,~E_1(x)=x-\frac{1}{2},~E_2(x)=x^2-x,~E_3(x)=x^3-\frac{3}{2}x^2+\frac{1}{4},~E_4(x)=x^4-2x^3+x.$$
The following properties are well known in the literature and will be used later.
\begin{itemize}
\item [$(P_1)$]
For $n\geq2,$  we have
$
B_n=\begin{cases}
0 &\text{ if } n\text{ is odd},\\
+ve &\text{ if } n\equiv2\ (mod\ 4),\\
-ve &\text { if } n\equiv0\ (mod\ 4).
\end{cases}
$
    \item [$(P_2)$] $E_n(0)=\frac{2}{n+1}(1-2^{n+1})B_{n+1}.$
    \item [$(P_3)$] $E_n^{\prime}(x)=nE_{n-1}(x)$ and $E_{n}(1-x)=(-1)^{n}E_n(x)$.
\end{itemize}

\begin{lem}\label{eulerlemma}
The Euler polynomials satisfy the following properties.
\begin{itemize}
       \item [$(i)$] For all $0<x<\frac{1}{2},$  we have $E_n(x)=\begin{cases}
+ve &\text{ if } n\equiv0,3\ (mod\ 4),\\
-ve &\text{ if } n\equiv1,2\ (mod\ 4).\\
    \end{cases}$
    \item[$(ii)$]  For all $\frac{1}{2}<x<1,$ we have $E_n(x)=\begin{cases}
+ve &\text{ if } n\equiv0,1\ (mod\ 4),\\
-ve &\text{ if } n\equiv2,3\ (mod\ 4).\\
    \end{cases}$
\end{itemize}
\end{lem}
\begin{pf}
$(i)$ We prove the result by induction on $n.$ It is easy to show that the result is true for $n=0,1,2,3.$
Assume that the result is true for $n=m-1.$
Now if $m=4l+k,$ for some $l\in\mathbb{N}$ and $k\in\{0, 1, 2, 3\},$ then
$$E_{4l+k}^{\prime}(x)=(4l+k)E_{4l+k-1}(x)=
\begin{cases}
+ve \text{ if } k=0, 1,\\
-ve \text{ if } k=2, 3.
\end{cases}$$ Therefore $E_{4l+k}$ is strictly monotone in $(0, \frac{1}{2})$ which implies that
$$E_{4l+k}(x)\begin{cases}
>E_{4l}(0)=0&\text{ if } k=0,\\
<E_{4l+1}(\frac{1}{2})=0&\text{ if } k=1,\\
<E_{4l+2}(0)=0&\text{ if } k=2,\\
>E_{4l+3}(\frac{1}{2})=0&\text{ if } k=3,
\end{cases}$$
for all $0<x<\frac{1}{2}.$ 
Hence the proof by induction.

 A proof of $(ii)$ follows from $(i)$ and  the second equality in  $(P_3)$.  
\end{pf}

For $s\in\mathbb{C},$ $b\in\mathbb{C}\setminus{\mathbb{Z}},$ and $\lambda\in\mathbb{R},$
the Lerch zeta function $L(\lambda, b, s)$ is defined by
$$L(\lambda, b, s)=\sum\limits_{n=0}^{\infty}\frac{e^{2\pi i\lambda n}}{(n+b)^s}.$$
This series converges absolutely for $Re(s)>1.$ The author in \cite{Fer} proved that for $s\in\mathbb{N},$ $\lambda\in [0, 1),$ and $b\in\mathbb{C}\setminus{\mathbb{Z}},$ 
\begin{equation}\label{lerch}
L(\lambda, b, s)+(-1)^sL(-\lambda,-b,s)=\frac{1}{b^s}+\frac{(-1)^{s-1}}{(s-1)!}\frac{d^{s-1}}{db^{s-1}}\left(\frac{\pi(\cot{\pi b}+i)}{e^{2\pi i\lambda b}}\right).
\end{equation}

Let us introduce the doubly infinite Lerch zeta function for $\lambda\in[0,1),$ $x\in\mathbb{C}\setminus{\mathbb{Z}},$  and $m>1$ as follows:
\begin{equation}\label{Htmx}
H(\lambda,x,m):=\sum\limits_{n\in\mathbb{Z}}\dfrac{e^{2\pi i\lambda n}}{(n-x)^{m}}.   
\end{equation}
By the definition of Lerch zeta function and \eqref{lerch}, we have
\begin{eqnarray}\label{htxmlerch}
H(\lambda,x,m)
&=&\sum\limits_{n=0}^{\infty}\dfrac{e^{2\pi i\lambda n}}{(n-x)^{m}}+(-1)^m\sum\limits_{n=1}^{\infty}\dfrac{e^{-2\pi i\lambda n}}{(n+x)^{m}}\nonumber\\
&=&\sum\limits_{n=0}^{\infty}\dfrac{e^{2\pi i\lambda n}}{(n-x)^{m}}+(-1)^m\sum\limits_{n=0}^{\infty}\dfrac{e^{-2\pi i\lambda n}}{(n+x)^{m}}-\dfrac{1}{(-x)^m}\nonumber\\
&=&L(\lambda,-x,m)+(-1)^mL(-\lambda,x,m)-\dfrac{1}{(-x)^m}\nonumber\\
&=&\dfrac{\pi}{(m-1)!}\frac{d^{m-1}}{dx^{m-1}}\left(\frac{-\cot{\pi x}+i}{e^{-2\pi i\lambda x}}\right).
\end{eqnarray}
\begin{lem}
Let us define 
\begin{equation}\label{p(j+1)}
P_0(x)=x \text{ and }
P_{j+1}(x)=-\left[(j+1)xP_j(x)+(1-x^2)P^{\prime}_j(x)\right],~j\geq0.
\end{equation}
Then
\begin{equation}\label{pjcos}
P_j(\cos{\pi x})=\pi^{-j}\sin^{j+1}{\pi x}\dfrac{d^j}{dx^j}\cot{\pi x}.
\end{equation}
\end{lem}
\begin{pf}
We prove the result by induction on $j$. The result is obvious for $j=0, 1$.
We now assume that \eqref{pjcos} is true for $j=k$. Then 
$$P_k^{\prime}(\cos{\pi x})=-\pi^{-k}(k+1)\cos{
\pi x}\sin^{k-1}{\pi x}\dfrac{d^k}{dx^k}\cot{\pi x}-\pi^{-k-1}\sin^{k}{\pi x}\dfrac{d^{k+1}}{dx^{k+1}}\cot{\pi x}.$$
From \eqref{p(j+1)}, we get 
\begin{align*}
&P_{k+1}(\cos{\pi x})\\
\hspace{0.5cm}
&=-\left[(k+1)\cos{\pi x}P_k(\cos{\pi x})+\sin^2{\pi x}P^{\prime}_k(\cos{\pi x})\right]\\
&=-(k+1)\cos{\pi x}\pi^{-k}\sin^{k+1}{\pi x}\dfrac{d^k}{dx^k}\cot{\pi x}\\
&-\sin^2{\pi x}\left(-\pi^{-k}(k+1)\cos{
\pi x}\sin^{k-1}{\pi x}\dfrac{d^k}{dx^k}\cot{\pi x}-\pi^{-k-1}\sin^{k}{\pi x}\dfrac{d^{k+1}}{dx^{k+1}}\cot{\pi x}\right)\\
&=\pi^{-(k+1)}\sin^{k+2}{\pi x}\dfrac{d^{k+1}}{dx^{k+1}}\cot{\pi x}.
\end{align*}
Hence \eqref{pjcos} is valid for all $j\geq1.$
\end{pf}

Let us define 
\begin{eqnarray}\label{fmx}
F_m(e^{i\pi x},\beta):=e^{m\pi ix}e^{-2\pi i\beta x}\frac{\sin^{m}{\pi x}}{\pi^{m-1}}\frac{d^{m-1}}{dx^{m-1}}\left(\left(i-\cot{\pi x}\right)e^{2\pi i\beta x}\right),
\end{eqnarray}
for all $\beta\in\mathbb{R}$, $x\in(0, 1),$ and $m\geq1$. 
By Leibniz rule, we have
\begin{align}\label{fmx2}
F_m(e^{i\pi x},\beta)&=e^{m\pi ix}\sin^{m}{\pi x}\sum\limits_{j=0}^{m-1}\binom{m-1}{j}(2i\beta)^{m-1-j}\pi^{-j}\frac{d^j}{dx^j}\left(i-\cot{\pi x}\right).
\end{align}
\begin{lem}\label{lemmafmi}
$F_m(i, \beta)
=2^{m-1}i^me^{m\pi i/2}E_{m-1}(\beta),$ for all $\beta\in\mathbb{R}.$
\end{lem}

\begin{pf}
It is proved in \cite{Kolbig} that
\begin{equation}\label{cot}
\frac{d^j}{dx^j}(\cot{\pi x})\Bigg|_{x=\frac{1}{2}}=-\frac{1}{\pi}(2\pi)^{j+1}(2^{j+1}-1)\frac{|B_{j+1}|}{j+1},
\end{equation}
where $B_n$ denotes the Bernoulli number.
When $x=\frac{1}{2},$ by \eqref{cot}, \eqref{fmx2} becomes  
\begin{eqnarray}\label{fmi}
F_m(i, \beta)&=&i^me^{m\pi i/2}(2\beta)^{m-1}-g_m(i, \beta),
\end{eqnarray}
where 
\begin{eqnarray}\label{gmi}
g_m(i, \beta)&=&-e^{m\pi i/2}(2i\beta)^{m}\sum\limits_{j=1}^{m-1}\binom{m-1}{j}(2i\beta)^{-j-1}\pi^{-j}\frac{1}{\pi}(2\pi)^{j+1}(2^{j+1}-1)\frac{|B_{j+1}|}{j+1}\nonumber\\
&=&-e^{m\pi i/2}(2i\beta)^{m}\sum\limits_{j=1}^{m-1}\binom{m-1}{j}\left(\frac{1}{i\beta}\right)^{j+1}(2^{j+1}-1)\frac{|B_{j+1}|}{j+1}\nonumber\\
&=&-e^{m\pi i/2}(2i\beta)^{m}\sum\limits_{j=1}^{m-1}\binom{m-1}{j}\left(\frac{-i}{\beta}\right)^{j+1}(2^{j+1}-1)\frac{|B_{j+1}|}{j+1}.
\end{eqnarray}

It is easy to check that $(-i)^{j+1}|B_{j+1}|=\begin{cases}
-B_{j+1} &\text{ if } j\equiv1, 3\ (mod\ 4),\\
0 &\text{ if } j\equiv0, 2\ (mod\ 4).
\end{cases}$

Thus we obtain from \eqref{gmi} that
\begin{eqnarray}\label{gmi1}
g_m(i, \beta)
&=&e^{m\pi i/2}(2i\beta)^{m}\sum\limits_{j=1}^{m-1}\binom{m-1}{j}\frac{1}{\beta^{j+1}}(2^{j+1}-1)\frac{B_{j+1}}{j+1}\nonumber\\
&=&e^{m\pi i/2}\frac{(2i\beta)^{m}}{m}\sum\limits_{j=1}^{m}\binom{m}{j+1}\frac{1}{\beta^{j+1}}(2^{j+1}-1)B_{j+1}\nonumber\\
&=&e^{m\pi i/2}\frac{(2i\beta)^{m}}{m}\sum\limits_{j=2}^{m}\binom{m}{j}\frac{1}{\beta^{j}}(2^{j}-1)B_{j}\nonumber\\
&=&e^{m\pi i/2}(2i\beta)^{m}\left\{\frac{1}{m}\sum\limits_{j=1}^{m}\binom{m}{j}\frac{1}{\beta^{j}}(2^{j}-1)B_{j}-\frac{B_1}{\beta}\right\}\nonumber\\
&=&e^{m\pi i/2}(2i\beta)^{m}\left\{\frac{-1}{2m}\sum\limits_{j=1}^{m}\binom{m}{j}\frac{1}{\beta^{j}}(2-2^{j+1})B_{j}+\frac{1}{2\beta}\right\}\nonumber\\
&=&e^{m\pi i/2}(2i\beta)^{m}\left\{\frac{-\beta^{-m}}{2m}\sum\limits_{j=1}^{m}\binom{m}{j}\beta^{m-j}(2-2^{j+1})B_{j}+\frac{1}{2\beta}\right\}\nonumber\\
&=&e^{m\pi i/2}(2i\beta)^{m}\left\{\frac{-\beta^{-m}}{2}E_{m-1}(\beta)+\frac{1}{2\beta}\right\}.
\end{eqnarray}
Substituting \eqref{gmi1} in \eqref{fmi}, we can get our result.
\end{pf}
\begin{thm}
The function $F_m(u,\beta)\equiv F_m(e^{i\pi x},\beta),$ $u=e^{i\pi x}$ satisfies the following recurrence relation: For $m\geq2$,  we have
\begin{equation}\label{recurrence}
F_m(u,\beta)=\left(\beta(u^2-1)-(m-1)u^2\right)F_{m-1}(u,\beta)+\frac{u(u^2-1)}{2}F_{m-1}^{\prime}(u,\beta)
\end{equation}
with $F_1(u,\beta)=-1.$
Here $F_m^{\prime}$ denotes the (partial) derivative of $F_{m}$ with respect to $u.$
\end{thm}
\begin{pf}
By \eqref{pjcos}, \eqref{fmx2} becomes
\begin{align}
F_m(e^{i\pi x},\beta)
&=(2i\sin{\pi x})^{m}\tfrac{\beta^{m-1}}{2}e^{m\pi ix}
-\sum\limits_{j=0}^{m-1}\tbinom{m-1}{j}(2i\beta\sin{\pi x})^{m-1-j}e^{m\pi ix}P_j(\cos{\pi x}).\nonumber
\end{align}
Taking $u=e^{i\pi x}$, we have
\begin{align}\label{fmtilde}
F_m(u,\beta)
&=\frac{\beta^{m-1}}{2}(u^2-1)^{m}-g_m(u,\beta),
\end{align}
where
\begin{equation}\label{gm1}
g_m(u,\beta)=\sum\limits_{j=0}^{m-1}\binom{m-1}{j}\beta^{m-1-j}u^{j+1}(u^2-1)^{m-1-j}P_j\left(\frac{u+u^{-1}}{2}\right). 
\end{equation}
Replacing $m$ by $m-1$ in \eqref{gm1}, we have
\begin{align}
g_{m-1}(u,\beta)&=\sum\limits_{j=0}^{m-2}\binom{m-2}{j}\beta^{m-2-j}u^{j+1}(u^2-1)^{m-2-j}P_j\left(\frac{u+u^{-1}}{2}\right)\label{g1}\\
&=\sum\limits_{j=0}^{m-2}\binom{m-1}{j}\frac{(m-1-j)}{(m-1)}\beta^{m-2-j}u^{j+1}(u^2-1)^{m-2-j}P_j\left(\frac{u+u^{-1}}{2}\right)\nonumber\\
&=\sum\limits_{j=0}^{m-1}\binom{m-1}{j}\beta^{m-2-j}u^{j+1}(u^2-1)^{m-2-j}P_j\left(\frac{u+u^{-1}}{2}\right)\nonumber\\
&-\frac{1}{(m-1)}\sum\limits_{j=0}^{m-1}\binom{m-1}{j}j\beta^{m-2-j}u^{j+1}(u^2-1)^{m-2-j}P_j\left(\frac{u+u^{-1}}{2}\right)\nonumber\\
&=\frac{g_m(u,\beta)}{\beta(u^2-1)}-\sum\limits_{j=1}^{m-1}\tfrac{(m-2)!}{(j-1)!(m-2-(j-1))!}\beta^{m-2-j}u^{j+1}(u^2-1)^{m-2-j}P_j\left(\tfrac{u+u^{-1}}{2}\right)\nonumber\\
&=\frac{g_m(u,\beta)}{\beta(u^2-1)}-\sum\limits_{j=0}^{m-2}\binom{m-2}{j}\beta^{m-3-j}u^{j+2}(u^2-1)^{m-3-j}P_{j+1}\left(\frac{u+u^{-1}}{2}\right).\nonumber\\\label{gmu}
\end{align}
Differentiating \eqref{g1} with respect to $u,$ we get
\begin{align}\label{gprime}
&g_{m-1}^{\prime}(u,\beta)\nonumber\\
&=\sum\limits_{j=0}^{m-2}\tbinom{m-2}{j}\beta^{m-2-j}\Big[\Big\{2(m-2-j)u^{j+2}(u^2-1)^{m-3-j}\nonumber\\
&+(j+1)u^j(u^2-1)^{m-2-j}\Big\}P_{j}\left(\tfrac{u+u^{-1}}{2}\right)+u^{j+1}(u^2-1)^{m-2-j}\left(\tfrac{1-u^{-2}}{2}\right)P_j^{\prime}\left(\tfrac{u+u^{-1}}{2}\right)\Big]\nonumber\\
&=2\sum\limits_{j=0}^{m-2}\tbinom{m-2}{j}\beta^{m-2-j}u^{j+1}(u^2-1)^{m-3-j}\Bigg[(m-1)uP_{j}\left(\tfrac{u+u^{-1}}{2}\right)\nonumber\\
&\hspace{4.5cm}-(j+1)\tfrac{u+u^{-1}}{2}P_{j}\left(\tfrac{u+u^{-1}}{2}\right)-\left(\tfrac{u-u^{-1}}{2i}\right)^2P_j^{\prime}\left(\tfrac{u+u^{-1}}{2}\right)\Bigg]\nonumber\\
&=2\sum\limits_{j=0}^{m-2}\tbinom{m-2}{j}\beta^{m-2-j}u^{j+1}(u^2-1)^{m-3-j}\Bigg[(m-1)uP_{j}\left(\tfrac{u+u^{-1}}{2}\right)+P_{j+1}\left(\tfrac{u+u^{-1}}{2}\right)\Bigg]\nonumber\\
&=\tfrac{2(m-1)u}{u^2-1}g_{m-1}(u,\beta)+\tfrac{2\beta}{u}\sum\limits_{j=0}^{m-2}\tbinom{m-2}{j}\beta^{m-3-j}u^{j+2}(u^2-1)^{m-3-j}P_{j+1}\left(\tfrac{u+u^{-1}}{2}\right),
\end{align}
using \eqref{p(j+1)} and \eqref{g1}.
From \eqref{gmu}, \eqref{gprime} becomes
\begin{equation*}
g_{m-1}^{\prime}(u,\beta)=\frac{2(m-1)u}{u^2-1}g_{m-1}(u,\beta)+\frac{2\beta}{u}\left(\frac{g_m(u,\beta)}{\beta(u^2-1)}-g_{m-1}(u,\beta)\right).
\end{equation*}
Rewriting the above equation, we get 
\begin{eqnarray}\label{gm}
g_m(u,\beta)=\left(\beta(u^2-1)-(m-1)u^2\right)g_{m-1}(u,\beta)+\frac{u(u^2-1)}{2}g_{m-1}^{\prime}(u,\beta).
\end{eqnarray}
By \eqref{fmtilde}, we have
\begin{equation}\label{gmu1}
g_m(u,\beta)
=\frac{\beta^{m-1}}{2}(u^2-1)^{m}-F_m(u,\beta).
\end{equation}
Replacing $m$ by $m-1$ in \eqref{gmu1}, we get
\begin{equation}\label{gm-1u}
g_{m-1}(u,\beta)
=\frac{\beta^{m-2}}{2}(u^2-1)^{m-1}-F_{m-1}(u,\beta).
\end{equation}
Differentiating the above equation with respect to $u,$ we get
\begin{equation}\label{gm-1uprime}
g_{m-1}^{\prime}(u,\beta)
=(m-1)\beta^{m-2}u(u^2-1)^{m-2}-F_{m-1}^{\prime}(u,\beta).
\end{equation}
Substituting \eqref{gmu1}, \eqref{gm-1u}, and \eqref{gm-1uprime} in \eqref{gm}, we get \eqref{recurrence}.
\end{pf}

By induction, we can easily show that the algebraic polynomial $F_m(u,\beta)$ defined in \eqref{recurrence} is an even function in a variable $u$. Let us define
\begin{equation}\label{gmfm}
G_m(t,\beta):=F_m(\sqrt{t},\beta).    
\end{equation}
Differentiating with respect to $t,$ we get
$$G_m^{\prime}(t,\beta)=\dfrac{1}{2\sqrt{t}}F_m^{\prime}(\sqrt{t},\beta).$$
Substituting $G_m$ and $G_m^\prime$  in \eqref{recurrence}, we get
\begin{equation}\label{recurrence1}
G_m(t,\beta)= [\beta(t-1)-(m-1)t]G_{m-1}(t,\beta)+t(t-1)G_{m-1}^{\prime}(t,\beta)   
\end{equation}
with $G_1(t,\beta)=-1.$
It is easy to check that $G_m(t,\beta)$ is a polynomial of degree $m-1$ in a variable $t$ and $G_{m}(1,\beta)=(-1)^m(m-1)!.$
The first few polynomials of $G_m(t, \beta)$ are
$$G_2(t, \beta)=(1-\beta)t+\beta,~G_3(t, \beta)=-(1-\beta)^2t^2-(1+2\beta-2\beta^2)t-\beta^2,$$
$$G_4(t, \beta)=(1-\beta)^3t^3+(3\beta^3-6\beta^2+4)t^2+(1+3\beta+3\beta^2-3\beta^3)t+\beta^3.$$

\begin{thm}\label{phigm}
The exponential splines $\Phi_{m-1}(\beta,t)$ satisfy the following properties.
\begin{itemize}
\item [$(i)$] $\Phi_{m-1}(\beta, t)=\dfrac{(-1)^m}{(m-1)!}G_m\left(\dfrac{1}{t}, \beta\right)$.
\item [$(ii)$] $t(m-1)\Phi_{m-1}(\beta, t)=\left(m-1+\beta(t-1)\right)\Phi_{m-2}(\beta, t)+t(1-t)\Phi_{m-2}^{\prime}(\beta, t).$
\end{itemize}
\end{thm}
\begin{pf}
From the definition of the Zak transform, \eqref{zak}, \eqref{htxmlerch}, and \eqref{fmx}, we have
\begin{eqnarray*}
\Phi_{m-1}(\beta,e^{-2\pi ix})
&=&\mathcal{Z}Q_m(\beta,-x)\\
&=&e^{-2\pi i \beta x}\sum\limits_{k\in\mathbb{Z}}\widehat{Q_{m}}(-x-k)e^{-2\pi i\beta k}\\
&=&e^{-2\pi i \beta x}\sum\limits_{k\in\mathbb{Z}}\sinc^m{(-k-x)}e^{-\pi im(-k-x)}e^{-2\pi i \beta k}\\
&=&e^{-2\pi i \beta x}e^{\pi i m x} \dfrac{\sin^m{\pi x}}{\pi^m}\sum\limits_{k\in\mathbb{Z}} \dfrac{e^{-2\pi i \beta k}}{(x+k)^m}\\
&=&(-1)^me^{-2\pi i \beta x}e^{\pi i m x} \dfrac{\sin^m{\pi x}}{\pi^m} H(\beta,x,m)\\
&=&(-1)^me^{-2\pi i \beta x}e^{\pi i m x} \dfrac{\sin^m{\pi x}}{(m-1)!\pi^{m-1}}\frac{d^{m-1}}{dx^{m-1}}\left(\frac{-\cot{\pi x}+i}{e^{-2\pi i\beta x}}\right)\\
&=&\dfrac{(-1)^m}{(m-1)!}F_m(e^{\pi i x}, \beta)\\
&=&\dfrac{(-1)^m}{(m-1)!}G_m(e^{2\pi i x}, \beta),
\end{eqnarray*}
which implies that
\begin{eqnarray}\label{phi1}
\Phi_{m-1}(\beta, t)=\frac{(-1)^m}{(m-1)!}G_m\left(\frac{1}{t},\beta\right).
\end{eqnarray}
Differentiating the above equation with respect to $t,$ we get
\begin{eqnarray}\label{phi2}
\Phi_{m-1}^{\prime}(\beta, t)=-\frac{(-1)^m}{(m-1)!}t^{-2}G_m^{\prime}\left(\frac{1}{t},\beta\right).
\end{eqnarray}
Replacing $t$ by $\frac{1}{t}$ in \eqref{recurrence1}, we get
\begin{eqnarray}
G_m\left(\frac{1}{t},\beta\right)
&=&\frac{1}{t}[\beta(1-t)-(m-1)]G_{m-1}\left(\frac{1}{t},\beta\right)+\frac{1}{t^2}(1-t)G_{m-1}^{\prime}\left(\frac{1}{t},\beta\right).\nonumber
\end{eqnarray}
Substituting \eqref{phi1} and \eqref{phi2} in the the above equation, we get
\begin{eqnarray*}
(-1)^m(m-1)!\Phi_{m-1}(\beta, t)&=&\frac{1}{t}[\beta(1-t)-(m-1)](-1)^{m-1}(m-2)!\Phi_{m-2}(\beta, t)\\
&&-\frac{1}{t^2}(1-t)(-1)^{m-1}(m-2)!t^2\Phi_{m-1}^{\prime}(\beta, t),
\end{eqnarray*}
from which we get the recurrence relation $(ii)$.
\end{pf}
\begin{lem}\label{recurrenceprop}
The polynomials  $G_m(t,\beta)$ satisfy the following properties:
\begin{itemize}
\item [$(i)$] For all $\beta\in\mathbb{R},$ we have
\begin{equation}\label{inverse}
t^{m-1}G_m\left(\frac{1}{t},1-\beta\right)=G_m(t,\beta).
\end{equation}
\item[$(ii)$]  For all $\beta>0,$ we have  $G_m(0,\beta)= \begin{cases}
+ve&\text{ if } m \text{ is even},\\
-ve&\text{ if } m \text{ is odd}.
\end{cases}$
\item[$(iii)$]  For all $0<\beta<\frac{1}{2},$ we have $G_m(-1,\beta)= \begin{cases}
-ve&\text{ if } m \equiv1,2~\ (mod\ 4),\\
+ve&\text{ if } m \equiv0,3~\ (mod\ 4).
\end{cases}$ 
\item[$(iv)$]  For all $1/2<\beta<1,$ we have  $G_m(-1,\beta)= \begin{cases}
-ve&\text{ if } m \equiv0, 1~\ (mod\ 4),\\
+ve&\text{ if } m \equiv2, 3~\ (mod\ 4).
\end{cases}$ 
\end{itemize}
\end{lem}

\begin{pf}
$(i)$
It is clear that \eqref{inverse} is true for $m=2$ because 
$$tG_2\left(\frac{1}{t},1-\beta\right)=t\left(\frac{1}{t}+\beta\left(1-\frac{1}{t}\right)\right)=G_2(t,\beta).$$
Now assume that $t^{m-1}G_m\left(\frac{1}{t},1-\beta\right)=G_m(t,\beta).$ Differentiating this equation with respect to $t,$ we have
\begin{eqnarray}
(m-1)t^{m-2}G_m\left(\frac{1}{t}, 1-\beta\right)-t^{m-3}G_m^{\prime}\left(\frac{1}{t}, 1-\beta\right)
=G_m^{\prime}(t,\beta),\nonumber
\end{eqnarray}
which implies that
\begin{eqnarray}\label{sign2}
t^{m-3}G_m^{\prime}\left(\dfrac{1}{t}, 1-\beta\right)
&=&(m-1)t^{m-2}G_m\left(\dfrac{1}{t}, 1-\beta\right)-G_m^{\prime}(t,\beta)\nonumber\\
&=&(m-1)t^{-1}G_m(t,\beta)-G_m^{\prime}(t,\beta).
\end{eqnarray}
Replacing $t$ by $\frac{1}{t},$ $\beta$ by $1-\beta,$ and $m$ by $m+1$ in \eqref{recurrence1},  we obtain from \eqref{sign2} that
\begin{eqnarray*}
&&\hspace{-1.5cm}t^mG_{m+1}\left(\frac{1}{t},1-\beta\right)\\
\hspace{1cm}&=&[(1-\beta)(1-t)-m]t^{m-1}G_m\left(\frac{1}{t},1-\beta\right)+t(1-t)t^{m-3}G_m^{\prime}\left(\frac{1}{t},1-\beta\right)\\
%&=&[(1-\beta)(1-t)-m]G_m(t,\beta)+t(1-t)(m-1)t^{-1}G_m(t,\beta)-t(1-t)G_m^{\prime}(t,\beta)\\
&=&[(1-\beta)(1-t)-m+(1-t)(m-1)]G_m(t,\beta)+t(t-1)G_m^{\prime}(t,\beta)\\
&=&[\beta(t-1)-mt]G_m(t,\beta)+t(t-1)G_m^{\prime}(t,\beta)=G_{m+1}(t,\beta).
\end{eqnarray*}
Hence, by induction \eqref{inverse} is proved for all $m.$

$(ii)$ Since $G_1(t,\beta)=-1,$ $G_1(0,\beta)=-ve.$
For $m=2,$ $G_2(t,\beta)=t+\beta(1-t).$
So, $G_2(0,\beta)=\beta>0,$ for all $\beta>0.$
Assume that $$G_m(0,\beta)=\begin{cases}
+ve&\text{ if } m \text{ is even},\\
-ve&\text{ if } m \text{ is odd}.
\end{cases}$$
From the recurrence relation  \eqref{recurrence1}, we have
$$G_{m+1}(t,\beta)=[\beta(t-1)-mt]G_m(t,\beta)+t(t-1)G_{m}^{\prime}(t,\beta),$$
and hence
\begin{eqnarray*}
 G_{m+1}(0,\beta)=-\beta G_m(0,\beta)=\begin{cases}
-ve&\text{ if } m+1 \text{ is odd},\\
+ve&\text{ if } m+1 \text{ is even}.
\end{cases}
\end{eqnarray*}
Hence the proof by induction.

Since $G_m(-1, \beta)=F_m(i, \beta),$ a proof of $(iii)$ and $(iv)$ follows from Lemma \ref{eulerlemma} and \ref{lemmafmi}.
\end{pf}

\begin{thm}\label{lemma2} For $m\geq2,$  $G_m(t,\beta)$ satisfies the following properties.
\begin{itemize}
\item[$(i)$] All the zeros of $G_m(t,\beta)$ are simple and negative for $0<\beta<1.$ 
\item[$(ii)$] If $0<\beta<\frac{1}{2},$ then 
$G_m(t,\beta)$ has exactly $\left\lfloor\dfrac{m}{2}\right\rfloor$ zeros in $|t|<1$ and $\left\lfloor\dfrac{m-1}{2}\right\rfloor$ zeros in $|t|>1.$
\item[$(iii)$] 
If $\frac{1}{2}<\beta<1,$  then 
$G_m(t,\beta)$ has exactly $\left\lfloor\dfrac{m-1}{2}\right\rfloor$ zeros in $|t|<1$ and $\left\lfloor\dfrac{m}{2}\right\rfloor$ zeros in $|t|>1.$
\end{itemize}
\end{thm}
\begin{pf}
We prove the result by induction on $m$.

$(i)$ Assume that $0<\beta<1$.  For $m=2,$ $G_2(t,\beta)=(1-\beta)t+\beta$ has only one zero which is negative.
%in $|t|<1$ and has no zeros in $|t|>1.$
Assume that $G_m(t,\beta)$ has $m-1$ zeros which are simple and negative. Let $\lambda_1, \lambda_2,\dots, \lambda_{m-1}$ be zeros of  $G_m(t,\beta)$ such that 
$$\lambda_{m-1}<\lambda_{m-2}<\dots<\lambda_1<0.$$
Here, the zeros depend on $\beta$; for simplicity, we omit 
$\beta$ in our notation.
By Lemma \ref{recurrenceprop}$(ii)$,  $G_m(t,\beta)$ is of the form
\begin{equation}\label{sign3}
G_m(t,\beta)=(-1)^mC_{m,\beta}(t-\lambda_1)\dots(t-\lambda_{m-1}), 
\end{equation}
for some positive constant $C_{m,\beta}$.
By \eqref{recurrence1}, we have
$$G_{m+1}(\lambda_\nu, \beta)=(\lambda_\nu-1)\lambda_\nu G_m^{\prime}(\lambda_\nu, \beta),~\nu=1,\dots, m-1.$$ 
From $\eqref{sign3}$, it is clear that
$(-1)^{\nu-1}(-1)^mG_m^{\prime}(\lambda_\nu,\beta)>0,$ for all $1\leq\nu\leq m-1.$
Hence
\begin{equation}\label{sign}
(-1)^{\nu-1}(-1)^mG_{m+1}(\lambda_\nu, \beta)>0,~\nu=1,\dots, m-1.  
\end{equation}
By Lemma \ref{recurrenceprop}$(ii)$, we have
\begin{equation}\label{sign1}
G_{m+1}(0,\beta)=\begin{cases}
-ve&\text{ if } m \text{ is even},\\
+ve&\text{ if } m \text{ is odd}.
\end{cases}
\end{equation}
We can conclude from Lemma \ref{recurrenceprop}$(i)$ that the coefficient of $t^m$ in $G_{m+1}(t,\beta)$ is $G_{m+1}(0,1-\beta).$ Therefore, by \eqref{sign1}, we have
$$\text{the coefficient of $t^m$ in $G_{m+1}(t,\beta)$}
=\begin{cases}
-ve&\text{ if } m \text{ is even},\\
+ve&\text{ if } m \text{ is odd}.
\end{cases}$$
Consequently, there always exists some $\mu<\lambda_{m-1}$ such that $G_{m+1}(\mu,\beta)<0.$ 
By intermediate value theorem, we can conclude that $G_{m+1}(t,\beta)$ has zeros in each of the $m$ intervals 
$$(-\infty, \lambda_{m-1}),(\lambda_{m-1}, \lambda_{m-2}),\dots,(\lambda_2, \lambda_1),(\lambda_1, 0).$$
$(ii)$ Assume that $G_{m}(t,\beta)$ has $\left\lfloor \dfrac{m}{2}\right\rfloor$ zeros in $|t|<1$ and $\left\lfloor \dfrac{m-1}{2}\right\rfloor$ zeros in $|t|>1$ for $0<\beta<\frac{1}{2}.$ Then
$$\lambda_{m-1}<\lambda_{m-2}<\dots<\lambda_{\left\lfloor \tfrac{m}{2}\right\rfloor+1}<-1<\lambda_{\left\lfloor \tfrac{m}{2}\right\rfloor}<\dots<\lambda_1<0.$$
\noindent
\underline{Case $1$: $m$ is odd}.  In this case, we have
$$G_{m+1}(\lambda_{\left\lfloor \tfrac{m}{2}\right\rfloor},\beta)= \begin{cases}
+ve&\text{ if } m \equiv1~\ (mod\ 4),\\
-ve&\text{ if } m \equiv3~\ (mod\ 4).
\end{cases} $$
By Lemma \ref{recurrenceprop}$(iii)$,
$G_{m+1}(t,\beta)$ has one zero in the interval $\left(-1, \lambda_{\left\lfloor \tfrac{m}{2}\right\rfloor}\right).$ Hence  $G_{m+1}(t,\beta)$ has $\left\lfloor \dfrac{m}{2}\right\rfloor+1=\left\lfloor\dfrac{m+1}{2}\right\rfloor$ zeros in $|t|<1$ 
and $m-\left\lfloor\dfrac{m+1}{2}\right\rfloor=\dfrac{m-1}{2}=\left\lfloor\dfrac{m}{2}\right\rfloor$ zeros in $|t|>1.$\\
\noindent
\underline{Case $2$: $m$ is even}. In this case, we have $$G_{m+1}(\lambda_{\left\lfloor \tfrac{m}{2}\right\rfloor},\beta)=\begin{cases}
-ve&\text{ if } m \equiv0~\ (mod\ 4),\\
+ve&\text{ if } m \equiv2~\ (mod\ 4).
\end{cases}$$
By Lemma \ref{recurrenceprop}$(iii)$, $G_{m+1}(t,\beta)$ has one zero in the interval $(\lambda_{\left\lfloor \tfrac{m}{2}\right\rfloor+1}, -1).$ Hence  $G_{m+1}(t,\beta)$ has $\left\lfloor \dfrac{m}{2}\right\rfloor=\left\lfloor\dfrac{m+1}{2}\right\rfloor$ zeros in $|t|<1$ and $m-\left\lfloor\dfrac{m}{2}\right\rfloor=\left\lfloor\dfrac{m}{2}\right\rfloor$ zeros in $|t|>1.$

A proof of $(iii)$ follows from $(ii)$ and Lemma \ref{recurrenceprop} $(i)$.
\end{pf}
\begin{cor}\label{zerooflerch}
$H(\lambda, x, m)\neq0,$  for all $\lambda\in(0, 1)\setminus\{\frac{1}{2}\},$ $x\in\mathbb{R}\setminus{\mathbb{Z}},$ and $m\geq2.$
\end{cor}
\begin{pf}
By \eqref{htxmlerch}, \eqref{fmx}, and \eqref{gmfm}, we have
\begin{eqnarray}
H(\lambda, x, m)=\frac{(2\pi i)^m}{(m-1)!}\frac{e^{2\pi i\lambda x}}{(e^{2\pi ix}-1)^m}G_m(e^{2\pi ix},\lambda).\nonumber
\end{eqnarray}
Hence we can conclude the result from Theorem \ref{lemma2}.
\end{pf}
\begin{cor}\label{zeroofphi}
For $m\geq2,$ $\Phi_{m-1}(\beta,t)$  satisfies the following properties.
\begin{itemize}
\item [$(i)$] All the zeros of $\Phi_{m-1}(\beta,t)$ are simple and negative for $0<\beta<1.$ 
\item[$(ii)$] For all $0<\beta<\frac{1}{2},$ $\Phi_{m-1}(\beta, t)$ has exactly $\lfloor\frac{m}{2}\rfloor$ zeros in $|t|>1$ and $\lfloor\frac{m-1}{2}\rfloor$ zeros in $|t|<1$.
\item[$(iii)$] For all $\frac{1}{2}<\beta<1,$ $\Phi_{m-1}(\beta, t)$ has exactly $\lfloor\frac{m-1}{2}\rfloor$ zeros in $|t|>1$ and $\lfloor\frac{m}{2}\rfloor$ zeros in $|t|<1$. 
\end{itemize}
\end{cor}
\begin{pf}
The result follows from Theorem \ref{phigm} and \ref{lemma2}.
\end{pf}
 
\section{Construction of complete interpolation sets for shift-invariant spline spaces}

The Euler-Frobenius polynomials $\Pi_{m}(t)$ and the modified Euler-Frobenius polynomials $\widetilde{\Pi}_{m}(t)$ 
are respectively  defined by
\begin{equation*}
\Pi_{m}(t)=m!\displaystyle\sum_{j=0}^{m-1}Q_{m+1}\left(j+1\right)t^j 
\text{ and }
\widetilde{\Pi}_{m}(t)=2^mm!\displaystyle\sum_{j=0}^{m}Q_{m+1}\left(j+\dfrac{1}{2}\right)t^j.
\end{equation*}
The polynomials $\Pi_{m}(t)$ and $\widetilde{\Pi}_{m}(t)$ are self-inversive monic polynomials and their zeros are simple and negative. 
They satisfy the following recurrence relations
\begin{equation}\label{euler}
\Pi_{m+1}(t)=(1+mt)\Pi_{m}(t)+t(1-t)\Pi_{m}^{\prime}(t),~ ~~\Pi_0(t)=1.
\end{equation}
and 
\begin{equation}\label{eulermodified}
\widetilde{\Pi}_{m+1}(t)=(1+(2m+1)t)\widetilde{\Pi}_{m}(t)+2t(1-t)\widetilde{\Pi}_{m}^{\prime}(t),~~\widetilde{\Pi}_{0}(t)=1.
\end{equation}
For more information about Euler-Frobenius and modified Euler-Frobenius polynomials, we refer to \cite{Schoenberg}.
\begin{lem}\label{lemma5.1}
$G_m(t,0)=(-1)^mt\Pi_{m-1}(t) \text{ and }G_m(t,\tfrac{1}{2})=\frac{(-1)^{m}}{2^{m-1}}\widetilde{\Pi}_{m-1}(t),$
where $\Pi_{m-1}(t)$ and $\widetilde{\Pi}_{m-1}(t)$ are defined in $\eqref{euler}$ and \eqref{eulermodified} respectively.   
\end{lem}
\begin{pf}    
For $m=2,$ $G_2(t,0)=t=(-1)^2t\Pi_1(t).$
Assume that $G_m(t,0)=(-1)^mt\Pi_{m-1}(t).$
Differentiating with respect to $t,$ we have 
$G_m^{\prime}(t,0)=(-1)^m(\Pi_{m-1}(t)+t\Pi_{m-1}^{\prime}(t)).$
By \eqref{recurrence1},
\begin{eqnarray*}
G_{m+1}(t,0)&=&-mtG_m(t,0)+t(t-1)G_m^{\prime}(t,0)\\
&=&(-1)^{m+1}t\{(mt-t+1)\Pi_{m-1}(t)+t(1-t)\Pi_{m-1}^{\prime}(t)\}\\
&=&(-1)^{m+1}t\Pi_{m}(t).
\end{eqnarray*}
Hence the result is true by induction.

Similarly, we can show that the second equality is valid for all $m\geq2$ by induction.
\end{pf}
\begin{thm}\label{cisspline}
Let $m\geq2$.
The set $\Lambda_{\alpha}=\left\{\left\langle\frac{m}{2}\right\rangle+n:n\in\mathbb{N}_0\}\cup\{\alpha+\left\langle\frac{m}{2}\right\rangle+n:n\in\mathbb{N}^{-}\right\}$ is a complete interpolation set for $V(Q_m)$ if and only if $|\alpha|<1/2.$ 
\end{thm}
\begin{pf}
By \eqref{symbol} and the definition of Zak transform, we have
\begin{eqnarray}\label{psiphigm}
\Psi_{\alpha,m}(x)
%&=&e^{-2\pi i(\alpha+a)x}\Psi_m^{\dagger}(x)^{-1}\sum\limits_{n\in\mathbb{Z}}\widehat{Q_m}(-x+n)e^{-2\pi i(-\alpha_0-a)n}\nonumber\\
&=&e^{-2\pi i(\alpha+a)x}\Psi^{\dagger}(x)^{-1}\mathcal{Z}\widehat{Q_m}(-x, -(\alpha_0+a))\nonumber\\
&=&e^{-2\pi i(\alpha+a) x}\Psi_{m}^\dagger(x)^{-1} \sum\limits_{n\in\mathbb{Z}}\sinc^m{(n-x)}e^{-\pi im(n-x)}e^{2\pi i(\alpha+a) n}.\nonumber
\end{eqnarray}
It is clear that $\Psi_{\alpha,m}(0)=\Psi_{\alpha,m}(1)=\Psi_{m}^\dagger(0)^{-1}\neq0.$
Let $\alpha_0=\langle\alpha\rangle$ and $\alpha_1=\langle \alpha_0+a\rangle.$ For $x\in(0, 1),$
 we obtain from \eqref{zak} and Theorem \ref{phigm} $(i)$ that
\begin{eqnarray}
\Psi_{\alpha,m}(x)&=&e^{-2\pi i(\alpha+a)x}\Psi^{\dagger}(x)^{-1}\mathcal{Z}Q_m(\alpha_0+a, -x)e^{2\pi i(\alpha_0+a)x}\nonumber\\
&=&e^{-2\pi i\lfloor\alpha\rfloor x}\Psi^{\dagger}(x)^{-1}\sum\limits_{n\in\mathbb{Z}}Q_m(\alpha_0+a-n)e^{-2\pi inx}\nonumber\\
&=&e^{-2\pi i\lfloor\alpha\rfloor x}e^{-2\pi i\lfloor\langle\alpha\rangle+a\rfloor x}\Psi^{\dagger}(x)^{-1}\sum\limits_{n\in\mathbb{Z}}Q_m(\alpha_1-n)e^{-2\pi inx}\nonumber\\
&=&e^{-2\pi i(\alpha+a-\alpha_1) x}\Psi_m^{\dagger}(x)^{-1}\Phi_{m-1}(\alpha_1,e^{-2\pi ix}),\nonumber\\
&=&(-1)^m\frac{\Psi_m^{\dagger}(x)^{-1}}{(m-1)!}\frac{1}{e^{2\pi i(\alpha-\alpha_1+a) x}}G_{m}(e^{2\pi ix},\alpha_1),
\end{eqnarray}

Since $G_{m}(1,\alpha_1)=(-1)^m(m-1)!,$ $\Psi_{\alpha,m}(x)$ is a continuous function on $[0,1]$.
Notice that
$$\Psi_m^{\dagger}(x)=\begin{cases}
\dfrac{e^{2\pi ix}\Pi_{m-1}(e^{2\pi ix})}{(m-1)!} &\text{ if } m \text{ is even},\\
\dfrac{\widetilde{\Pi}_{m-1}(e^{2\pi ix})}{2^{m-1}(m-1)!} &\text{ if } m \text{ is odd}.
\end{cases}$$
Let $\Psi_{\alpha,m}(x)\equiv\varphi_{\alpha,m}(z)$, where $z=e^{2\pi ix}$. Then \eqref{psiphigm} becomes
\begin{align}\label{eqn1}
\Psi_{\alpha,m}(x)=&\varphi_{\alpha,m}(z)
=\begin{cases}
\dfrac{G_{m}(z,\alpha_0)}{z^{\lfloor\alpha\rfloor+1}\Pi_{m-1}(z)}&\text{ if }  m \text{ is even},\\\\
-\dfrac{2^{m-1}G_{m}(z,\alpha_1)}{z^{\alpha-\alpha_1+\frac{1}{2}}\widetilde{\Pi}_{m-1}(z)} &\text{ if }  m \text{ is odd}.
\end{cases}
\end{align}

By Lemma \ref{recurrenceprop} $(i)$ and Lemma \ref{lemma5.1}, it is easy to show that  $\Psi_{\alpha, m}(1/2)=0$ if $\langle\alpha\rangle=\frac{1}{2}$  So, we can conclude from Theorem \ref{cis} that $\Lambda_{\alpha}$ is not a complete interpolation set if $\langle\alpha\rangle=\frac{1}{2}.$

If $\langle\alpha\rangle\in(0, 1)\setminus\{1/2\},$ then 
it follows from Theorem \ref{lemma2} $(ii)$ and $(iii)$ that $G_m(z,\alpha_0)\neq0$ and $G_m(z,\alpha_1)\neq0$ for all $z$ in the unit circle $\mathbb{T}$. 
When $\langle\alpha\rangle=0,$ by Lemma \ref{lemma5.1},  $G_m(z,0)\neq0$ for all $z\in\mathbb{T}$ if $m$ is even and  $G_m(z,\tfrac{1}{2})\neq0$ for all $z\in\mathbb{T}$ if $m$ is odd.
Hence $\Psi_{\alpha,m}(x)\neq0$  for all $x\in [0,1]$ if  $\langle\alpha\rangle\neq\frac{1}{2}.$

Since the function $w=\varphi_{\alpha,m}(z)$ maps $\mathbb{T}$ onto the closed curve $\Psi_{\alpha,m}$ in the $w$- plane,
by the argument principle from complex analysis \cite{Ahlfors},
we find that
\begin{eqnarray}\label{eqn2}
\mathrm{Ind}(\Psi_{\alpha,m},0)=\dfrac{1}{2\pi i}\int\limits_{\Psi_{\alpha,m}}\dfrac{1}{w}dw
=\dfrac{1}{2\pi i}\int\limits_{|z|=1}\dfrac{\varphi_{\alpha,m}'(z)}{\varphi_{\alpha,m}(z)}dz=Z_{\alpha,m}-P_{\alpha,m},
\end{eqnarray}
where $Z_{\alpha,m}$ and $P_{\alpha,m}$ denote respectively the number of zeros and poles of $\varphi_{\alpha,m}(z)$ inside the unit circle, counting multiplicities.\\
\noindent
\underline{Case $1$. $\langle\alpha\rangle=0$}.
By Lemma \ref{lemma5.1}, \eqref{eqn1} becomes
\begin{eqnarray*}
\varphi_{\alpha,m}(z)
=\dfrac{1}{z^{\lfloor\alpha\rfloor}},
\end{eqnarray*}
which implies that $\mathrm{Ind}(\Psi_{\alpha,m},0)=0$ if and only if $\alpha=0$.\\ 
\underline{Case $2$.  $\langle\alpha\rangle\in(0, 1)\setminus\{\frac{1}{2}\}$ and $m$ is even}. 
By Theorem \ref{lemma2}, $G_{m}(z,\langle\alpha\rangle)$ has $\frac{m}{2}$ zeros in $|z|<1$ for all $0<\langle\alpha\rangle<\frac{1}{2}$ and $\frac{m-2}{2}$ zeros in $|z|<1$ for all $\frac{1}{2}<\langle\alpha\rangle<1.$
Since $\Pi_{m-1}(z)$ has exactly $\frac{m-2}{2}$ zeros in $|z|<1$,
it follows from \eqref{eqn1} and \eqref{eqn2}  that
\begin{eqnarray*}
\mathrm{Ind}(\Psi_{\alpha,m},0)  
=\begin{cases}
-\lfloor\alpha\rfloor&\text{ if }0<\langle\alpha\rangle<\frac{1}{2},\\
-1-\lfloor\alpha\rfloor&\text{ if }\frac{1}{2}<\langle\alpha\rangle<1.
\end{cases}\\
\end{eqnarray*}
Now it is easy to conclude that $\mathrm{Ind}(\Psi_{\alpha,m},0)=0$ if and only if $0<|\alpha|<\frac{1}{2}.$\\
\underline{Case $3$. $\langle\alpha\rangle\in(0, 1)\setminus\{\frac{1}{2}\}$ and $m$ is odd  }. 
If $\alpha_0=\langle \alpha\rangle\in\left(0, \frac{1}{2}\right)$, then $\alpha_1=\alpha_0+\frac{1}{2}\in \left(\frac{1}{2},1\right)$
and if 
$\alpha_0\in\left(\frac{1}{2}, 1\right)$, then $\alpha_1=\alpha_0-\frac{1}{2}\in\left(0, \frac{1}{2}\right).$
By Theorem \ref{lemma2},
$G_{m}(z,\alpha_1)$ has $\frac{m-1}{2}$ zeros in $|z|<1$ for all $\langle\alpha\rangle\in(0,1)\setminus\{\frac{1}{2}\}$. Since  $\widetilde{\Pi}_{m-1}(z)$ has exactly $\frac{m-1}{2}$ zeros in $|z|<1,$ it follows from \eqref{eqn1} and \eqref{eqn2} that
\begin{eqnarray*}
\mathrm{Ind}(\Psi_{\alpha,m},0)
=\begin{cases}
-\lfloor\alpha\rfloor&\text{ if } 0<\langle\alpha\rangle<\frac{1}{2},\\
 -1 -\lfloor\alpha\rfloor&\text{ if }\frac{1}{2}<\langle\alpha\rangle<1,
\end{cases} 
\end{eqnarray*}
and hence $\mathrm{Ind}(\Psi_{\alpha,m},0)$ is zero if and only if $0<|\alpha|<\frac{1}{2}$.

Now we can conclude our result from Theorem $\ref{cis}.$
\end{pf}
\section*{Conclusion}
This paper focuses on characterizing the complete interpolation property of the set $(a+\mathbb{N}_0)\cup(\alpha+a+\mathbb{N}^{-})$ for shift-invariant spaces using  Toeplitz operators. 
We characterize all $\alpha$ for which $\mathbb{N}_0\cup\alpha+\mathbb{N}^{-}$ forms a complete interpolation set for a transversal-invariant space.
The theory we develop introduces a new type of recurrence relation for exponential splines, explores the location of zeros of these splines, and identifies the zero-free region of the doubly infinite Lerch zeta function. These findings are then utilized to prove that $\left\langle\frac{m}{2}\right\rangle+\mathbb{N}_0\cup\alpha+\left\langle\frac{m}{2}\right\rangle+\mathbb{N}^{-}$ is a complete interpolation set for a  shift-invariant spline spaces of order $m\geq 2$ if and only if $|\alpha|<1/2$. 
The results and ideas of this paper may be applied in number theory for understanding the zeros of special functions, in approximation theory to enhance methods for function approximation and smoothing, and in time-frequency analysis for Gabor frames.

\end{document}